%VERSION: 21 November 2005 (without date, for arXiv)
\magnification\magstep1
%		*****	  MSSYMB.TeX	*****		       8 Jul 87
%
%	This file contains the definitions for the symbols in the two
%	"extra symbols" fonts created at the American Math. Society.

\catcode`\@=11
\def\relaxnext@{\let\next\relax}
\font\tenmsx=msam10       %changed from msxm to msam by H.S.
\font\sevenmsx=msam7      %changed from msxm to msam by H.S.
\font\fivemsx=msam5       %changed from msxm to msam by H.S.
\font\tenmsy=msbm10       %changed from msym to msbm by H.S.
\font\sevenmsy=msbm7      %changed from msym to msbm by H.S.
\font\fivemsy=msbm5       %changed from msym to msbm by H.S.
\newfam\msxfam
\newfam\msyfam
\textfont\msxfam=\tenmsx  \scriptfont\msxfam=\sevenmsx
  \scriptscriptfont\msxfam=\fivemsx
\textfont\msyfam=\tenmsy  \scriptfont\msyfam=\sevenmsy
  \scriptscriptfont\msyfam=\fivemsy

\def\hexnumber@#1{\ifcase#1 0\or1\or2\or3\or4\or5\or6\or7\or8\or9\or
	A\or B\or C\or D\or E\or F\fi }

%  The following 13 lines establish the use of the Euler Fraktur font.
%  To use this font, remove % from beginning of these lines.
\font\teneuf=eufm10
\font\seveneuf=eufm7
\font\fiveeuf=eufm5
\newfam\euffam
\textfont\euffam=\teneuf
\scriptfont\euffam=\seveneuf
\scriptscriptfont\euffam=\fiveeuf
\def\frak{\relaxnext@\ifmmode\let\next\frak@\else
 \def\next{\errmessage{Use \string\frak\space only in math
mode}}\fi\next}
\def\goth{\relaxnext@\ifmmode\let\next\frak@\else 
\def\next{\errmessage{Use \string\goth\space only in math
mode}}\fi\next} 
\def\frak@#1{{\frak@@{#1}}}
\def\frak@@#1{\fam\euffam#1}
%  End definition of Euler Fraktur font.

\edef\msx@{\hexnumber@\msxfam}
\edef\msy@{\hexnumber@\msyfam}

\mathchardef\boxdot="2\msx@00
\mathchardef\boxplus="2\msx@01
\mathchardef\boxtimes="2\msx@02
\mathchardef\square="0\msx@03
\mathchardef\blacksquare="0\msx@04
\mathchardef\centerdot="2\msx@05
\mathchardef\lozenge="0\msx@06
\mathchardef\blacklozenge="0\msx@07
\mathchardef\circlearrowright="3\msx@08
\mathchardef\circlearrowleft="3\msx@09
\mathchardef\rightleftharpoons="3\msx@0A
\mathchardef\leftrightharpoons="3\msx@0B
\mathchardef\boxminus="2\msx@0C
\mathchardef\Vdash="3\msx@0D
\mathchardef\Vvdash="3\msx@0E
\mathchardef\vDash="3\msx@0F
\mathchardef\twoheadrightarrow="3\msx@10
\mathchardef\twoheadleftarrow="3\msx@11
\mathchardef\leftleftarrows="3\msx@12
\mathchardef\rightrightarrows="3\msx@13
\mathchardef\upuparrows="3\msx@14
\mathchardef\downdownarrows="3\msx@15
\mathchardef\upharpoonright="3\msx@16

\mathchardef\downharpoonright="3\msx@17
\mathchardef\upharpoonleft="3\msx@18
\mathchardef\downharpoonleft="3\msx@19
\mathchardef\rightarrowtail="3\msx@1A
\mathchardef\leftarrowtail="3\msx@1B
\mathchardef\leftrightarrows="3\msx@1C
\mathchardef\rightleftarrows="3\msx@1D
\mathchardef\Lsh="3\msx@1E
\mathchardef\Rsh="3\msx@1F
\mathchardef\rightsquigarrow="3\msx@20
\mathchardef\leftrightsquigarrow="3\msx@21
\mathchardef\looparrowleft="3\msx@22
\mathchardef\looparrowright="3\msx@23
\mathchardef\circeq="3\msx@24
\mathchardef\succsim="3\msx@25
\mathchardef\gtrsim="3\msx@26
\mathchardef\gtrapprox="3\msx@27
\mathchardef\multimap="3\msx@28
\mathchardef\therefore="3\msx@29
\mathchardef\because="3\msx@2A
\mathchardef\doteqdot="3\msx@2B

\mathchardef\triangleq="3\msx@2C
\mathchardef\precsim="3\msx@2D
\mathchardef\lesssim="3\msx@2E
\mathchardef\lessapprox="3\msx@2F
\mathchardef\eqslantless="3\msx@30
\mathchardef\eqslantgtr="3\msx@31
\mathchardef\curlyeqprec="3\msx@32
\mathchardef\curlyeqsucc="3\msx@33
\mathchardef\preccurlyeq="3\msx@34
\mathchardef\leqq="3\msx@35
\mathchardef\leqslant="3\msx@36
\mathchardef\lessgtr="3\msx@37
\mathchardef\backprime="0\msx@38
\mathchardef\risingdotseq="3\msx@3A
\mathchardef\fallingdotseq="3\msx@3B
\mathchardef\succcurlyeq="3\msx@3C
\mathchardef\geqq="3\msx@3D
\mathchardef\geqslant="3\msx@3E
\mathchardef\gtrless="3\msx@3F
\mathchardef\sqsubset="3\msx@40
\mathchardef\sqsupset="3\msx@41
\mathchardef\vartriangleright="3\msx@42
\mathchardef\vartriangleleft="3\msx@43
\mathchardef\trianglerighteq="3\msx@44
\mathchardef\trianglelefteq="3\msx@45
\mathchardef\bigstar="0\msx@46
\mathchardef\between="3\msx@47
\mathchardef\blacktriangledown="0\msx@48
\mathchardef\blacktriangleright="3\msx@49
\mathchardef\blacktriangleleft="3\msx@4A
\mathchardef\vartriangle="0\msx@4D
\mathchardef\blacktriangle="0\msx@4E
\mathchardef\triangledown="0\msx@4F
\mathchardef\eqcirc="3\msx@50
\mathchardef\lesseqgtr="3\msx@51
\mathchardef\gtreqless="3\msx@52
\mathchardef\lesseqqgtr="3\msx@53
\mathchardef\gtreqqless="3\msx@54
\mathchardef\Rrightarrow="3\msx@56
\mathchardef\Lleftarrow="3\msx@57
\mathchardef\veebar="2\msx@59
\mathchardef\barwedge="2\msx@5A
\mathchardef\doublebarwedge="2\msx@5B
\mathchardef\angle="0\msx@5C
\mathchardef\measuredangle="0\msx@5D
\mathchardef\sphericalangle="0\msx@5E
\mathchardef\varpropto="3\msx@5F
\mathchardef\smallsmile="3\msx@60
\mathchardef\smallfrown="3\msx@61
\mathchardef\Subset="3\msx@62
\mathchardef\Supset="3\msx@63
\mathchardef\Cup="2\msx@64

\mathchardef\Cap="2\msx@65

\mathchardef\curlywedge="2\msx@66
\mathchardef\curlyvee="2\msx@67
\mathchardef\leftthreetimes="2\msx@68
\mathchardef\rightthreetimes="2\msx@69
\mathchardef\subseteqq="3\msx@6A
\mathchardef\supseteqq="3\msx@6B
\mathchardef\bumpeq="3\msx@6C
\mathchardef\Bumpeq="3\msx@6D
\mathchardef\lll="3\msx@6E

\mathchardef\ggg="3\msx@6F

\mathchardef\circledS="0\msx@73
\mathchardef\pitchfork="3\msx@74
\mathchardef\dotplus="2\msx@75
\mathchardef\backsim="3\msx@76
\mathchardef\backsimeq="3\msx@77
\mathchardef\complement="0\msx@7B
\mathchardef\intercal="2\msx@7C
\mathchardef\circledcirc="2\msx@7D
\mathchardef\circledast="2\msx@7E
\mathchardef\circleddash="2\msx@7F
\def\ulcorner{\delimiter"4\msx@70\msx@70 }
\def\urcorner{\delimiter"5\msx@71\msx@71 }
\def\llcorner{\delimiter"4\msx@78\msx@78 }
\def\lrcorner{\delimiter"5\msx@79\msx@79 }
\def\yen{\mathhexbox\msx@55 }
\def\checkmark{\mathhexbox\msx@58 }
\def\circledR{\mathhexbox\msx@72 }
\def\maltese{\mathhexbox\msx@7A }
\mathchardef\lvertneqq="3\msy@00
\mathchardef\gvertneqq="3\msy@01
\mathchardef\nleq="3\msy@02
\mathchardef\ngeq="3\msy@03
\mathchardef\nless="3\msy@04
\mathchardef\ngtr="3\msy@05
\mathchardef\nprec="3\msy@06
\mathchardef\nsucc="3\msy@07
\mathchardef\lneqq="3\msy@08
\mathchardef\gneqq="3\msy@09
\mathchardef\nleqslant="3\msy@0A
\mathchardef\ngeqslant="3\msy@0B
\mathchardef\lneq="3\msy@0C
\mathchardef\gneq="3\msy@0D
\mathchardef\npreceq="3\msy@0E
\mathchardef\nsucceq="3\msy@0F
\mathchardef\precnsim="3\msy@10
\mathchardef\succnsim="3\msy@11
\mathchardef\lnsim="3\msy@12
\mathchardef\gnsim="3\msy@13
\mathchardef\nleqq="3\msy@14
\mathchardef\ngeqq="3\msy@15
\mathchardef\precneqq="3\msy@16
\mathchardef\succneqq="3\msy@17
\mathchardef\precnapprox="3\msy@18
\mathchardef\succnapprox="3\msy@19
\mathchardef\lnapprox="3\msy@1A
\mathchardef\gnapprox="3\msy@1B
\mathchardef\nsim="3\msy@1C
%\mathchardef\napprox="3\msy@1D
\mathchardef\ncong="3\msy@1D

\mathchardef\varsubsetneq="3\msy@20
\mathchardef\varsupsetneq="3\msy@21
\mathchardef\nsubseteqq="3\msy@22
\mathchardef\nsupseteqq="3\msy@23
\mathchardef\subsetneqq="3\msy@24
\mathchardef\supsetneqq="3\msy@25
\mathchardef\varsubsetneqq="3\msy@26
\mathchardef\varsupsetneqq="3\msy@27
\mathchardef\subsetneq="3\msy@28
\mathchardef\supsetneq="3\msy@29
\mathchardef\nsubseteq="3\msy@2A
\mathchardef\nsupseteq="3\msy@2B
\mathchardef\nparallel="3\msy@2C
\mathchardef\nmid="3\msy@2D
\mathchardef\nshortmid="3\msy@2E
\mathchardef\nshortparallel="3\msy@2F
\mathchardef\nvdash="3\msy@30
\mathchardef\nVdash="3\msy@31
\mathchardef\nvDash="3\msy@32
\mathchardef\nVDash="3\msy@33
\mathchardef\ntrianglerighteq="3\msy@34
\mathchardef\ntrianglelefteq="3\msy@35
\mathchardef\ntriangleleft="3\msy@36
\mathchardef\ntriangleright="3\msy@37
\mathchardef\nleftarrow="3\msy@38
\mathchardef\nrightarrow="3\msy@39
\mathchardef\nLeftarrow="3\msy@3A
\mathchardef\nRightarrow="3\msy@3B
\mathchardef\nLeftrightarrow="3\msy@3C
\mathchardef\nleftrightarrow="3\msy@3D
\mathchardef\divideontimes="2\msy@3E
\mathchardef\varnothing="0\msy@3F
\mathchardef\nexists="0\msy@40
\mathchardef\mho="0\msy@66
\mathchardef\eth="0\msy@67
\mathchardef\eqsim="3\msy@68
\mathchardef\beth="0\msy@69
\mathchardef\gimel="0\msy@6A
\mathchardef\daleth="0\msy@6B
\mathchardef\lessdot="3\msy@6C
\mathchardef\gtrdot="3\msy@6D
\mathchardef\ltimes="2\msy@6E
\mathchardef\rtimes="2\msy@6F
\mathchardef\shortmid="3\msy@70
\mathchardef\shortparallel="3\msy@71
\mathchardef\smallsetminus="2\msy@72
\mathchardef\thicksim="3\msy@73
\mathchardef\thickapprox="3\msy@74
\mathchardef\approxeq="3\msy@75
\mathchardef\succapprox="3\msy@76
\mathchardef\precapprox="3\msy@77
\mathchardef\curvearrowleft="3\msy@78
\mathchardef\curvearrowright="3\msy@79
\mathchardef\digamma="0\msy@7A
\mathchardef\varkappa="0\msy@7B
\mathchardef\hslash="0\msy@7D
\mathchardef\hbar="0\msy@7E
\mathchardef\backepsilon="3\msy@7F
% Use the next 4 lines with AMS-TeX:
%\def\Bbb{\relaxnext@\ifmmode\let\next\Bbb@\else
% \def\next{\err@{Use \string\Bbb\space only in math mode}}\fi\next}
%\def\Bbb@#1{{\Bbb@@{#1}}}
%\def\Bbb@@#1{\noaccents@\fam\msyfam#1}
% Use the next 4 lines if NOT using AMS-TeX:
\def\Bbb{\ifmmode\let\next\Bbb@\else
 \def\next{\errmessage{Use \string\Bbb\space only in math mode}}\fi\next}
\def\Bbb@#1{{\Bbb@@{#1}}}
\def\Bbb@@#1{\fam\msyfam#1}

\catcode`\@=12
%%%%%%%%%%%END OF MSSYMB
%%%%%%%%%%%lay-out macros:
\font\bigbf=cmbx12 scaled 1200
\def\Sect#1#2{\par\bigskip\goodbreak\line{\bf\hbox 
 to\parindent{#1\hfil}#2\hfil}}

\def\Satz#1#2{\medskip\noindent{\bf #1 #2.} \begingroup\it}
\def\endSatz{\endgroup\par\medskip}
\def\Thm#1{\Satz{Theorem}{#1}}
\def\endThm{\endSatz}
\def\Prop#1{\Satz{Proposition}{#1}}
\def\endProp{\endSatz}
\def\Lem#1{\Satz{Lemma}{#1}}
\def\endLem{\endSatz}
\def\Cor#1{\Satz{Corollary}{#1}}
\def\endCor{\endSatz}
\def\Def#1{\medskip\noindent{\bf Definition #1. }}
\def\endDef{\par\medskip}

\def\Proof{\smallskip\noindent{\bf Proof. }}
\def\Proofof#1{\smallskip\noindent{\bf Proof of #1. }}
\def\endProof{\ignorespaces\ $\square$\par\smallskip}

\def\Eq#1{\eqno(\hbox{\rm #1})}

%%%%%%math macros:
\def\complex{{\Bbb C}}
\def\real{{\Bbb R}}
\def\Natural{{\Bbb N}}
\def\integer{{\Bbb Z}}
\def\BbbD{{\Bbb D}}

%defines := with : as a relation

\def\cA{{\cal A}}
\def\cB{{\cal B}}
\def\cC{{\cal C}}
\def\cF{{\cal F}}
\def\cH{{\cal H}}
\def\cJ{{\cal J}}

\def\cL{{\cal L}}
\def\cM{{\cal M}}
\def\cN{{\cal N}}
\def\cO{{\cal O}}
\def\cP{{\cal P}}

\def\Fou{{\cF}}
\def\Hyp{{\cH}}
\def\Pal{{\cP}}

\def\Mer{{\cM}}
\def\Par{{\cP}}

\def\Wave{{\cJ}}

\def\Cusp{{\cA}}%{{\cC}}
\def\Lau{{\cL}}

\def\supp{\mathop{\rm supp}}

\def\Hom{\mathop{\rm Hom}\nolimits}

\def\Re{\mathop{\rm Re}}

\def\laur{{\mathop{\rm laur}}}

\def\AC{\mathop{\rm AC}\nolimits}

\def\PW{\mathop{\rm PW}\nolimits}
%{\mathop{\rm H}}

\def\iq{{\rm q}}

\def\inp#1#2{\langle#1,#2\rangle}
\def\frac#1#2{{{#1}\over{#2}}}
\def\ay{{\fa}}

\def\Cci{C_c^\infty}
\def\Cinf{C^\infty}

\newdimen\prestarheight\prestarheight=4pt
\newdimen\scriptprestarheight\scriptprestarheight=2.8pt
\newdimen\scriptscriptprestarheight\scriptscriptprestarheight=2pt
\def\prestar#1#2{\mathchoice
 {\raise#1\prestarheight\hbox{$\scriptstyle *$}\kern-.#2em}
 {\raise#1\prestarheight\hbox{$\scriptstyle *$}\kern-.#2em}
 {\raise#1\scriptprestarheight\hbox{$\scriptscriptstyle *$}\kern-.#2em}
 {\raise#1\scriptscriptprestarheight\hbox{$\scriptscriptstyle *$}\kern-.#2em}}

\def\col{\colon}
\def\CciXt{\Cci(X\col\tau)}

\def\faqdF{\fa_{F\iq}^*}

\def\nE{E^\circ}

\def\DGH{{\BbbD}(G/H)}

\def\Aq{A_\iq}
\def\faq{\fa_\iq}

\def\faqd{\faq^*}
\def\faqcd{\fa_{\iq\siC}^*}

\def\oC{{}^\circ{\cC}}
\def\East{E^*}

\def\ayF{{\ay_{F\iq}}}

\def\faqF{\ayF}

\def\fgc{\fg_{\siC}}

\def\baraqd{\bar{\fa}^*_\iq}

%{{\rm E}}

%%%%%%extra by Erik:

\def\iC{{\scriptstyle\complex}}
\def\siC{{\scriptscriptstyle\complex}}

\def\siC{{\scriptstyle\complex}}
\def\faqdc{\fa_{\iq\siC}^*}
\def\faqFdperpc{\fa_{F\iq \siC}^{*\perp}}

\def\faqFdc{\fa_{F\iq \siC}^*}
\def\faqFd{\fa_{F\iq}^*}

%%%nye til PW
\def\Laufu{\Mer(\faqcd,\gS)^*_\laur}
\def\LauoC{\Laufu\otimes\oC^*}

\def\ACR{\AC_\real(X\col\tau)}

\def\CF{\Cusp_F}

%^{\stt}}

\def\PWXt{\PW(X\col\tau)}

\def\PWMX{\PW_{\! M}(X\col\tau)}

%%%%%%%%new for distPW
\def\Cmi{C^{-\infty}}
\def\Ccmi{\Cmi_c}
\def\distc{\Ccmi(X\col\tau)}
\def\distM{\Cmi_M(X\col\tau)}
\def\fbk{\fb_{\rm k}}
\def\fbc{\fb_{\iC}}
\def\Wb{{W(\fb)}}
\def\PWMd{\PW_M^*(X\col\tau)}
\def\PWd{\PW^*(X\col\tau)}

\def\hpl{Y}
%%%%greek alphabet
\def\ga{\alpha}
\def\gd{\delta}
\def\geps{\epsilon}

\def\gl{\lambda}
\def\gL{\Lambda}
\def\gf{\varphi}

\def\gS{\Sigma}
%%%%fraktur
\def\fa{\frak a}
\def\fb{\frak b}
\def\fg{\frak g}
\def\fh{\frak h}

\def\fk{\frak k}

\def\fp{\frak p}
\def\fq{\frak q}
%%%%LABELS:
%sections:
\newcount\sectcount\sectcount=1
\def\sectdef{\number\sectcount}\def\nxsect{\advance\sectcount by 1}
\edef\SectZ{\sectdef}\nxsect
\edef\SectA{\sectdef}\nxsect
\edef\SectY{\sectdef}\nxsect
\edef\SectB{\sectdef}\nxsect
\edef\SectM{\sectdef}\nxsect
\edef\SectC{\sectdef}\nxsect
\edef\SectD{\sectdef}\nxsect
\edef\SectE{\sectdef}\nxsect
\edef\SectL{\sectdef}\nxsect
\edef\SectG{\sectdef}\nxsect
\edef\SectJ{\sectdef}\nxsect
\edef\SectF{\sectdef}\nxsect
\edef\SectH{\sectdef}\nxsect
\edef\SectI{\sectdef}\nxsect
\edef\SectK{\sectdef}\nxsect
%lemmas, propositions etc:
\newcount\propcount\propcount=1
\def\propdef{\number\propcount}\def\nxprop{\advance\propcount by 1}
\propcount=1
\edef\DefBA{\SectB.\propdef}\nxprop
\edef\LemBA{\SectB.\propdef}\nxprop
\edef\DefBC{\SectB.\propdef}\nxprop
\edef\DefBD{\SectB.\propdef}\nxprop
\edef\DefBE{\SectB.\propdef}\nxprop
\edef\ThmBF{\SectB.\propdef}\nxprop
\propcount=1
\edef\LemID{\SectM.\propdef}\nxprop
\propcount=1
\edef\LemCA{\SectC.\propdef}\nxprop
\propcount=1
\edef\ThmEA{\SectE.\propdef}\nxprop
\edef\LemEC{\SectE.\propdef}\nxprop
\edef\LemEB{\SectE.\propdef}\nxprop
\propcount=1
\edef\PropFA{\SectF.\propdef}\nxprop
\propcount=1
\edef\LemGA{\SectL.\propdef}\nxprop
\edef\LemHA{\SectL.\propdef}\nxprop
\edef\LemJA{\SectL.\propdef}\nxprop
\propcount=1
\edef\LemExtra{\SectJ.\propdef}\nxprop
\edef\LemJB{\SectJ.\propdef}\nxprop
\propcount=1
\edef\LemHB{\SectH.\propdef}\nxprop
\edef\CorHA{\SectH.\propdef}\nxprop
\propcount=1
\edef\ThmKA{\SectI.\propdef}\nxprop
\edef\LemIA{\SectI.\propdef}\nxprop
\propcount=1
\edef\DefIA{\SectK.\propdef}\nxprop
\edef\LemIB{\SectK.\propdef}\nxprop
\edef\PropIC{\SectK.\propdef}\nxprop
%\propcount=1
%\edef\ThmKB{\SectK.\propdef}\nxprop
%\edef\PropJB{\SectJ.\propdef}\nxprop
%equations:
\newcount\tagcount\tagcount=1
\def\tagdef{\number\tagcount}\def\nxtag{\advance\tagcount by 1}
\edef\Eqac{\SectA.\tagdef}\nxtag
\edef\Eqaa{\SectA.\tagdef}\nxtag
\edef\Eqab{\SectA.\tagdef}\nxtag
\tagcount=1\tagcount=1
\edef\Eqya{\SectY.\tagdef}\nxtag
\edef\Eqyb{\SectY.\tagdef}\nxtag
\tagcount=1\tagcount=1
\edef\Eqba{\SectB.\tagdef}\nxtag
\edef\Eqbd{\SectB.\tagdef}\nxtag
\nxtag
\nxtag
\tagcount=1\tagcount=1
\edef\Eqma{\SectM.\tagdef}\nxtag
\edef\Eqmb{\SectM.\tagdef}\nxtag
\tagcount=1\tagcount=1
\edef\Eqcb{\SectC.\tagdef}\nxtag
\edef\Eqca{\SectC.\tagdef}\nxtag
\edef\Eqcc{\SectC.\tagdef}\nxtag
\tagcount=1\tagcount=1
\edef\Eqdc{\SectD.\tagdef}\nxtag
\edef\Eqdd{\SectD.\tagdef}\nxtag
\edef\Eqda{\SectD.\tagdef}\nxtag
\edef\Eqdb{\SectD.\tagdef}\nxtag
\tagcount=1\tagcount=1
\edef\Eqea{\SectE.\tagdef}\nxtag
\edef\Eqeb{\SectE.\tagdef}\nxtag
\tagcount=1\tagcount=1
\edef\Eqgb{\SectL.\tagdef}\nxtag
\tagcount=1\tagcount=1
\edef\Eqgd{\SectG.\tagdef}\nxtag
\edef\Eqgc{\SectG.\tagdef}\nxtag
\edef\Eqga{\SectG.\tagdef}\nxtag
\tagcount=1\tagcount=1
\edef\Eqjc{\SectJ.\tagdef}\nxtag
\edef\Eqja{\SectJ.\tagdef}\nxtag
\tagcount=1\tagcount=1
\edef\Eqfa{\SectF.\tagdef}\nxtag
\edef\Eqfb{\SectF.\tagdef}\nxtag
\tagcount=1\tagcount=1
\edef\Eqhb{\SectH.\tagdef}\nxtag
\edef\Eqha{\SectH.\tagdef}\nxtag
\tagcount=1\tagcount=1
\edef\Eqia{\SectI.\tagdef}\nxtag
\edef\Eqkb{\SectI.\tagdef}\nxtag
\edef\Eqka{\SectI.\tagdef}\nxtag
\edef\Eqkc{\SectI.\tagdef}\nxtag
%References
\newcount\refcount
\def\refdef{\number\refcount}
\def\nxref{\advance\refcount by 1}
\refcount=1
\edef\rArth{\refdef}\nxref
\edef\rBsurvey{\refdef}\nxref
\edef\rBFJS{\refdef}\nxref
\edef\rBSmc{\refdef}\nxref
\edef\rBSrc{\refdef}\nxref
\edef\rBSFI{\refdef}\nxref
\edef\rBSaf{\refdef}\nxref
\edef\rBSPlI{\refdef}\nxref
\edef\rBSPlII{\refdef}\nxref
\edef\rBSPW{\refdef}\nxref
\edef\rBSvD{\refdef}\nxref
\edef\rDadok{\refdef}\nxref
\edef\rEHO{\refdef}\nxref
\edef\rGang{\refdef}\nxref
\edef\rHePW{\refdef}\nxref
\edef\rHeBAMS{\refdef}\nxref
\edef\rHorm{\refdef}\nxref
\edef\rSchaefer{\refdef}\nxref
\edef\rSbook{\refdef}\nxref
\edef\rSsurvey{\refdef}\nxref
%\edef\rWall{\refdef}\nxref

%%%%%%%%% END OF LAY OUT MACROS
\bigskip
\centerline{\bigbf A Paley-Wiener theorem for distributions}
\vskip5pt
\centerline{\bigbf on reductive symmetric spaces}
\bigskip
\centerline{by E. P. van den Ban and H. Schlichtkrull}
%\bigskip
%\centerline{21.11.2005}
\bigskip
\noindent{\bf Abstract}

Let $X=G/H$ be a reductive symmetric space and $K$ a maximal
compact subgroup of $G$. We study Fourier transforms of 
compactly supported $K$-finite distributions on $X$ and
characterize the image of the space of such distributions. 
\bigskip

\Sect\SectZ{Introduction}
The well-known Paley-Wiener-Schwartz theorem for the Fourier transform on
Euclidean space $\real^n$ characterizes the Fourier image of the
space $C_c^\infty(\real^n)$
of compactly supported smooth functions. The image is the space
of entire functions $\varphi\in\cO(\complex^n)$ with decay of
exponential type. The theorem has a counterpart, also 
well-known and also called the Paley-Wiener-Schwartz theorem,
where smooth functions are replaced by distributions, and where the
exponential decay condition is replaced by a similar 
exponential condition of
slow growth, see [\rHorm], Thm.\ 7.3.1.
The theorem for smooth functions was generalized to the 
Fourier transform of a reductive symmetric space $G/H$
in [\rBSPW]. It is the purpose of the present paper to establish
an analogue of the theorem for distributions in the same 
spirit and generality. 

In the more restricted case of a Riemannian symmetric space $G/K$,  
where $H=K$ is compact, a Paley-Wiener theorem for $K$-invariant
smooth functions was obtained by work of Helgason and Gangolli, 
[\rHePW], [\rGang],
and for general smooth functions by Helgason [\rHeBAMS].
A counterpart for distributions was given by Eguchi, Hashizume
and Okamoto in [\rEHO]. A different proof of the latter result
is given in [\rDadok].

Another important special case is that of a reductive Lie group,
considered as a symmetric space. In this case the Paley-Wiener
theorem of [\rBSPW] specializes to a theorem of Arthur [\rArth],
which describes the Fourier image of the space of compactly 
supported $K$-finite smooth functions on $G$ ($K$ being a maximal
compact subgroup). The theorem for distributions, which is
obtained in the present paper, is new in this `group case'.
The specialization to the group case is described in
[\rBSvD], to which we refer for further details.

The Paley-Wiener theorem of [\rBSPW] describes the Fourier
image of the space of $K$-finite compactly supported smooth
functions by an exponential type condition, combined with
a set of so-called Arthur-Campoli conditions. In the present 
theorem the exponential type condition is replaced by a 
condition of slow growth which is similar to its Euclidean analogue, 
whereas the additional Arthur-Campoli conditions 
remain the same as in [\rBSPW]. The precise statement of our main result
is given in Theorem \ThmBF, and its proof is given in
Sections \SectM-\SectH. 
The main tool in the proof is a
Fourier inversion formula, through which a function is
determined from its Fourier transform by means of certain
`residual' operators. Given a function $\varphi$ in the conjectured image 
space, we construct the distribution $f$, which is the candidate
for the inverse Fourier image, by means of this formula. 
The proof that $f$ has compact support and transforms to $\varphi$
is carried out by regularization with a Dirac sequence.
The mentioned inversion formula is generalized to
distributions in Corollary \CorHA, after the proof of Theorem \ThmBF.
Finally, in Sections \SectI-\SectK{}  
we discuss the topology on the
image space by which the Fourier transform becomes a
topological isomorphism.

For general background about harmonic analysis 
and Paley-Wiener theorems on reductive symmetric spaces
we refer to the survey articles [\rBsurvey], [\rSsurvey].

\medskip
\noindent{\it Acknowledgement.} We are grateful to 
Erik Thomas for helpful discussions related to the material
in Section  \SectI.

\Sect\SectA{Notation}
As in [\rBSPW] we use the notation and basic assumptions from [\rBSmc],
Sect.\ 2-3, 5-6 and [\rBSFI], Sect.\ 2. Only the most essential notions
will be recalled.

Let $G$ be a real reductive Lie group of Harish-Chandra's class,
and let $H$ be an open subgroup of the group of fixed points
for an involution $\sigma$. Then $X=G/H$ is a reductive symmetric space.
Let $K$ be a maximal compact subgroup of $G$, invariant under
$\sigma$, and let $\theta$ denote the corresponding Cartan involution.
Let $\fg$ denote the Lie algebra of $G$, which decomposes
in $\pm1$ eigenspaces for $\sigma$ and $\theta$ as $\fg=\fh+\fq=\fk+\fp$.
Then $\fh$ and $\fk$ are the Lie algebras of $H$ and $K$.
Let $\faq$ be a maximal abelian subspace of $\fq\cap\fp$, and choose
a positive system $\Sigma^+$ for the root system $\Sigma$ of
$\faq$ in $\fg$. This positive system determines a parabolic subgroup
$P$ of $G$, which will be fixed throughout the paper.
We also fix a finite dimensional unitary representation $(\tau,
V_\tau)$ of $K$. The normalized Eisenstein integrals associated with these
choices are denoted by $\nE(\psi\col\gl)\colon X\to V_\tau$, where
$\psi\in\oC$ and $\gl\in\faqdc$, as in [\rBSmc] p.\ 283.
Here $\oC=\oC(\tau)$ is the finite dimensional Hilbert space defined
in [\rBSmc], eq.\ (5.1).
The Eisenstein integrals depend
linearly on the parameter $\psi$ in this space, 
and as functions on $X$ they
belong to the space
$C^\infty(X\col\tau)$ of smooth $V_\tau$-valued functions on $X$ which are
$\tau$-spherical, that is, which satisfy the transformation rule
$$f(kx)=\tau(k)f(x),\qquad k\in K, x\in X.\Eq\Eqac$$
The 
adjoint of the linear map
$\psi\mapsto \nE(\psi\col-\bar\gl\col x)$ is denoted by $\East(\gl\col x)$,
see [\rBSFI], eq.\ (2.3), and the {\it Fourier transform}
for the $K$-type $\tau$ on $G/H$ is then defined by
$$
\Fou f(\gl)=\int_X \East(\gl\col x)f(x)\,dx \in\oC
\Eq\Eqaa$$
for $\gl\in\faqdc$ and for
$f$ in the space $\CciXt$ of compactly supported functions in
$C^\infty(X\col\tau)$, cf.\ [\rBSPW] Eq.~(2.1). Here $dx$ is an
invariant measure on $G/H$, normalized as in [\rBSmc] Section 3. 

The normalized Eisenstein integrals $\nE(\psi\col\gl\col x)$
depend meromorphically on the parameter $\gl\in\faqdc$, in a uniform way
with respect to the parameters $\psi$ and $x$. The nature of this
meromorphic dependence is crucial. It can be described as follows.
By a real $\gS$-configuration in $\faqdc$ we mean a
locally finite collection $\Hyp$ of affine
hyperplanes $\hpl$ in $\faqdc$ of the form $\hpl=\{\gl\mid \inp{\gl}{\ga_\hpl}=s_\hpl\}$,
where $\ga_\hpl\in\Sigma$ and $s_\hpl\in\real$.
Let $d\colon\Hyp\to\Natural$ be an arbitrary map. 
For $\omega\subset\faqdc$ we write 
$$\Hyp(\omega)=\{ Y\in\Hyp\mid Y\cap\omega\neq \emptyset\}$$
and, if the set $\Hyp(\omega)$ is finite, 
$$\pi_{\omega,d}(\gl)=\prod_{\hpl\in\Hyp(\omega)}
(\inp{\gl}{\ga_\hpl}-s_\hpl)^{d(\hpl)}.$$

Let $V$ be an arbitrary complete locally
convex vector space.
The linear space of meromorphic functions 
$\varphi\colon\faqdc\to V$, such that
$\pi_{\omega,d}\varphi$ is holomorphic on $\omega$ for all bounded
open sets $\omega\subset\faqdc$, is denoted
$\Mer(\faqd,\Hyp,d, V)$. 
It follows from [\rBSPW] Lemma 2.1 
and [\rBSFI] Prop.\ 3.1 that there exist 
a real $\Sigma$-configuration $\Hyp$ and a map $d\colon\Hyp\to\Natural$
such that the normalized Eisenstein integrals
$\gl\mapsto\nE(\psi\col\gl)$ belong to
$\Mer(\faqd,\Hyp,d, C^\infty(X)\otimes V_\tau)$
for all $\psi\in\oC$.

Clearly the dualized Eisenstein integrals $\East(\gl\col x)$ have
the same type of meromorphic dependence on $\gl$. 
We define 
$\Hyp=\Hyp(X,\tau)$ and $d=d_{X,\tau}\colon \Hyp\to\Natural$
as in [\rBSPW],
Section 2. Then $\Hyp$ is a real $\gS$-configuration, and
the map $\East\colon\gl\mapsto\East(\gl)$ satisfies
$$\East\in\Mer(\faqd,\Hyp,d, C^\infty(X)\otimes V_\tau^*\otimes\oC).
\Eq\Eqab$$
Moreover, $\Hyp(X,\tau)$ and $d_{X,\tau}$ are minimal with
respect to this property.

\Sect\SectY{The Fourier transform of a distribution}
The concept of `distributions' used in this paper is
that of {\it generalized functions}. By this we mean the
following. A generalized function on a smooth manifold
$X$ is a continuous linear form on the space of
compactly supported smooth densities on $X$.
We denote by $\Cmi(X)$ the space of generalized
functions on $X$, and by $\Ccmi(X)$ the subspace of generalized
functions with compact support. If a nowhere vanishing
smooth density $dx$ is given on
$X$, then the multiplication with $dx$ induces linear isomorphisms of
the spaces $\Cmi(X)$ and $\Ccmi(X)$ onto the topological
linear duals of $\Cci(X)$ and $C^\infty(X),$ 
respectively. If $f\in\Cmi(X)$ and $\phi\in\Cci(X)$, or if
$f \in \Ccmi(X)$ and $\phi \in\Cinf(X),$
then we write accordingly:
$$ \int_{X} \phi(x) \, f(x) \, dx = f\,dx(\phi).\Eq\Eqya
$$

Let $X=G/H$, as in Section 2, be  equipped with the invariant
measure $dx$. Then $dx$ is a nowhere vanishing smooth density. 
For $f\in\Cmi(X)$ the continuity of $f\,dx$, as a linear form on
$\Cci(X)$, can be expressed as follows. Let $X_1,\dots,X_n$ be
a linear basis for $\fg$, and for $\alpha=(\alpha_1,\dots,\alpha_n)$
a multi-index let $X^\alpha=X_1^{\alpha_1}\dots X_n^{\alpha_n}\in
\cal U(\fg)$. Then for each compact $\Omega\subset X$ there exist
constants $C,k$ such that 
$$\big|\int_X \phi(x)\,f(x)\,dx\big|\leq C\sup_{|\alpha|\leq k,x\in\Omega}
|L_{X^\alpha}\phi(x)|\Eq\Eqyb$$
for all $\phi\in\Cinf(X)$ with support in $\Omega$.

For each
positive number $M$ we denote by $\Cmi_M(X)$
the space of generalized functions with support in the compact
set $K\exp B_M H$. Here $B_M$ is the closed ball in $\faq$ 
centered at $0$ and of radius $M$. In view of 
the generalized Cartan 
decomposition $G=K\Aq H$ we have $\Ccmi(X)=\cup_M \Cmi_M(X)$.

A generalized function on $X$ with values in $V_\tau$ 
is called $\tau$-spherical
if it 
satisfies  
(\Eqac).
We denote by $\Cmi(X\col\tau)$
the space of $\tau$-spherical generalized functions on
$X$, and by $\distc$ and $\distM$ 
the subspaces of $\tau$-spherical
distributions with compact support, respectively
with support in 
$K\exp B_M H$.

If $f\in\Cmi(X\col\tau)$ and 
$\phi\in\Cci(X)\otimes V_\tau^*$, or if $f\in\distc$ and
$\phi\in\Cinf(X)\otimes  V_\tau^*$, 
then equation (\Eqya) still has a natural interpretation.
Via this pairing (\Eqya) we have thus established linear isomorphisms
of $\Cmi(X\col\tau)$ and $\distc$ with the topological linear
duals of $\Cci(X\col\tau^*)$ and $\Cinf(X\col\tau^*)$,
respectively.

Having established (\Eqya) in this generality
we can define 
the Fourier transform $\Fou f(\lambda)$ for $f\in\distc$
and $\gl\in\faqdc$
by the very same formula (\Eqaa)
by which it was defined for
$f\in\Cci(X\col\tau)$. The Fourier transform $\Fou f(\lambda)$  
is a $\oC$-valued meromorphic function
of $\gl$, and it follows from (\Eqab) that
$$\Fou f\in \Mer(\faqd,\Hyp(X,\tau),d_{X,\tau},\oC).$$

\Sect\SectB{The distributional Paley-Wiener space}

Recall the following definitions from [\rBSPW].

\Def\DefBA
Let $\Hyp$ be a real $\gS$-configuration in $\faqdc$, and
let $d\in\Natural^\Hyp$.
By $\Pal(\faqd,\Hyp,d)$ we denote the linear space of functions
$\gf\in\Mer(\faqd,\Hyp,d)$ with polynomial decay in the
imaginary directions, that is
$$\nu_{\omega,n}(\varphi):=
\sup_{\gl\in\omega+i\faqd} (1+|\gl|)^n\|\pi_{\omega,d}(\gl)\gf(\gl)\|
<\infty\Eq\Eqba$$
for all compact $\omega\subset\faqd$ and all $n\in\Natural$. 
The union of these spaces over all $d\colon\Hyp\to\Natural$
is denoted $\Pal(\faqd,\Hyp)$.

The space $\Pal(\faqd,\Hyp,d)$ is a Fr\' echet space with the topology 
defined by means of the seminorms $\nu_{\omega,n}$ in
(\Eqba).
\endDef

We recall the following result from
[\rBSPW], Lemma 3.7.

\Lem\LemBA Fourier
transform is continuous
$$\Fou\colon \CciXt\to\Pal(\faqd,\Hyp(X,\tau),d_{X,\tau})\otimes\oC.$$
\endLem

For $R\in\real$ we define
$\baraqd(P,R)=\{\gl\in\faqdc\mid
\forall\ga\in\gS^+:\Re\langle\gl,\ga\rangle<R\}$.

\Def\DefBC  Let $\Hyp=\Hyp(X,\tau)$ and $d=d_{X,\tau}$.
Let $\pi=\pi_{\baraqd(P,0),d}$. 
For each $M>0$ we define $\PWMX$ 
as the space of functions $\gf\in\Pal(\faqd,\Hyp,d)\otimes\oC$ 
for which 

\smallskip
(i) $\cL\gf=0$ 
for all $\cL\in\ACR$,

\smallskip
(ii)
$\sup_{\gl\in\baraqd(P,0)} (1+|\gl|)^n e^{-M\,|\!\Re\gl|}
\|\pi(\gl)\gf(\gl)\|<\infty$
for each $n\in\Natural$

\smallskip 
\noindent
(see [\rBSPW] Defn.\ 3.1 for the definition of 
$\ACR$).
Furthermore, the {\it Paley-Wiener space} $\PWXt$ is 
defined as
$$\PWXt=\cup_{M>0}\PWMX.$$
\endDef

The main result of [\rBSPW], Thm.\ 3.6, 
asserts that the Fourier transform 
is a linear isomorphism of $\Cinf_M(X\col\tau)$ onto 
$\PWMX$ for each $M>0$, and hence also of $\CciXt$
onto $\PWXt$.

We now introduce the following definitions. If $\Omega$ is a
topological space, we denote by $\cC(\Omega)$ the
collection of compact subsets $\omega\subset\Omega$, and by $\cN(\Omega)$
the set of maps $n\colon\cC(\Omega)\to\Natural$. 
%such that $\omega_1\subset\omega_2 \Rightarrow n(\omega_1)\leq n(\omega_2)$.

\Def\DefBD
Let $\Hyp$ be a real $\gS$-configuration in $\faqdc$, and
let $d\in\Natural^\Hyp$, $n\in\cN(\faqd)$.
By $\Pal^*(\faqd,\Hyp,d,n)$ we denote the linear space of functions
$\gf\in\Mer(\faqd,\Hyp,d)$ with at most polynomial growth of order
$n$ in the imaginary directions, that is
$$\nu^*_{\omega,n}(\varphi):=\sup_{\gl\in\omega+i\faqd} (1+|\gl|)^{-n(\omega)}
\|\pi_{\omega,d}(\gl)\gf(\gl)\|
<\infty\Eq\Eqbd$$
for all $\omega\in\cC(\faqd)$. 
The union $\cup_{n} \Pal^*(\faqd,\Hyp,d,n)$
is denoted by $\Pal^*(\faqd,\Hyp,d)$, and
the union $\cup_{d}\Pal^*(\faqd,\Hyp,d)$
is denoted by $\Pal^*(\faqd,\Hyp)$.

The space $\Pal^*(\faqd,\Hyp,d,n)$ is a 
locally convex topological vector %Fr\'echet 
space with the topology 
defined by means of the seminorms $\nu^*_{\omega,n}$ in
(\Eqbd).
This topological vector space is discussed further 
in Section \SectK{}, where it is shown to be
Fr\'echet under a natural condition on the map $n$.
However, this property is not needed at present. 
\endDef

\Def\DefBE Let $\Hyp=\Hyp(X,\tau)$ and $d=d_{X,\tau}$.
For each $M>0$ 
we define $\PWMd$ 
as the space of functions $\gf\in\Pal^*(\faqd,\Hyp,d)\otimes\oC$ 
for which 

\smallskip
(i) $\cL\gf=0$ 
for all $\cL\in\ACR$

\smallskip
(ii)
$\sup_{\gl\in\baraqd(P,0)} (1+|\gl|)^{-n} e^{-M\,|\!\Re\gl|}
\|\pi(\gl)\gf(\gl)\|<\infty$
for some $n\in\Natural$.

\smallskip 
\noindent
The {\it Paley-Wiener space} $\PWd$ is 
then defined by
$\PWd=\cup_{M>0}\PWMd.$
\endDef

It is clear that $\PWMX\subset\PWMd$ for all $M>0$ and 
that $\PWXt\subset\PWd$.

We can now state our main theorem.

\Thm\ThmBF The Fourier transform $\Fou$ is a linear isomorphism of\/
$\distM$ onto the Paley-Wiener space $\PWMd$ for each $M>0$, and hence also 
of\/ $\distc$ onto $\PWd$.
\endThm

The proof will be given in the course of the following 
Sections \SectM-\SectH.

\Sect\SectM{An estimate}
Let a real $\gS$-configuration $\Hyp$ in $\faqdc$
and a map $d\colon \Hyp\to\Natural$ be given.
Let $V$ be a finite dimensional normed vector space.

\Lem\LemID 
Let  $\omega_0\subset \omega_1\subset\faqd$, and assume that
$\Hyp(\omega_1)$ is finite.
Assume also that for some
$\delta>0$ the open set 
$$\omega=\{\gl+\mu\mid \gl\in\omega_1, |\mu|<\delta\},$$
is contained in $\omega_1$
and satisfies $\Hyp(\omega)=\Hyp(\omega_0)$ (for example, this
condition is fulfilled if 
$\omega_1$ is compact and contained in the interior 
of $\omega_2$). 

Let $p \in \Pi_{\gS}(\faqd)$, $n\in\integer$ and 
$M\geq 0$ be given.
There exists a constant $C >0$ such that 
$$
\eqalign{
\sup_{\omega_0 + i \faqd} &(1 + |\gl|)^{n} 
e^{-M|\!\Re(\gl)|}\|\pi_{\omega_0,d}\,(\gl)\gf(\gl)\|\cr
&\leq C 
\sup_{\omega_1 + i \faqd} (1 + |\gl|)^{n} 
e^{-M|\!\Re(\gl)|}\|\pi_{\omega_1,d}(\gl)\,p(\gl)\gf(\gl)\|}
\Eq\Eqma
$$
for all $\gf \in \Mer(\faqd, \Hyp, d, V)$.
\endLem

\Proof 
We first prove the result under the assumption that 
$d=0$ on $\Hyp(\omega_1)$. Then  $\pi_{\omega_0,d}=\pi_{\omega_1,d}=1$
and every function from $\Mer(\faqd, \Hyp, d, V)$
is holomorphic on the open set $\omega+i\faqd$.

It suffices to prove the
estimate for $p(\gl) = \inp{\ga}{\gl} - s,$ with $\ga \in \gS$ and
$s \in \complex.$ 
%We may assume $\delta<1$. 
We
fix $\mu \in \faqd$ such that $|\mu| < \gd$ and $\inp{\ga}{\mu} = c > 0.$
Then for every $\gl \in \omega_0 + i \faqd$ and $0 < r \leq 1,$
we have, by Cauchy's integral formula,
$$
%\eqalign{
\gf(\gl) 
%& 
=  \frac{1}{2\pi i} \int_{|z| = r} \frac{\gf(\gl + z \mu)}{z} \; dz
%\cr& 
=  \frac{1}{2\pi i} \int_{|z| = r} \frac{p(\gl + z \mu)
\gf(\gl + z \mu)}{(p(\gl) + cz)z} \; dz.
%}
$$
If $|p(\gl)| > 2c/3,$ we fix $r = 1/3 $, and if  $|p(\gl)| \leq  2c/3,$
we fix $r = 1 .$ In all cases we have $|p(\gl) + c z | \geq c/3$
for $|z|  = r.$ 
Using the above integral formula we thus obtain the estimate
$$
(1 + |\gl|)^n\|\gf(\gl)\| 
\leq \frac{3}{c}\,(1 + |\gl|)^n\,
\sup_{|z|=r} \|[p\gf](\gl + z \mu)\|
$$
We now observe that, for all $\gl \in \faqdc$ and 
$z \in \complex,\; |z| \leq 1,$
$$
(1-\delta)(1+|\gl|)\leq 1+|\gl|-\delta\leq 1+|\gl+z\mu|
$$
and hence
$$(1 + |\gl|)^n  \|\gf(\gl)\| 
\leq 
\frac{3}{c(1-\delta)} \,
\sup_{|z|=r} (1+|\gl+z\mu|)^n \|[p\gf](\gl + z \mu)\|.
$$
Since $|\!\Re(\gl+z\mu)|\leq|\!\Re\gl|+\delta$ for all $z$ with
$|z|=r$ we further obtain
$$\eqalign{
(1 + |\gl|)^n &e^{-M|\!\Re(\gl)|} \|\gf(\gl)\| \cr
\leq 
&\frac{3e^{M\delta}}{c(1-\delta)}  \,
\sup_{|z|=r} 
(1+|\gl+z\mu|)^n e^{-M|\!\Re(\gl+z\mu)|}\|[p\gf](\gl + z \mu)\|.}
$$
Now (\Eqma) follows, and we can
proceed to the general case.

From $\omega_0 \subset \omega_1$ it follows that
$\pi_{\omega_1, d} = q\pi_{\omega_0, d}$
with $q \in \Pi_\gS(\faqd).$ 
Define $d'\colon \Hyp\to\Natural$ by $d'=d$ on 
$\Hyp\setminus\Hyp(\omega_0)$ and $d'=0$ on $\Hyp(\omega_0)$.
By application of the first part of the proof, with $pq$, $\omega$ and
$d'$ in place of $p$, $\omega_1$ and $d$, 
there exists a constant $C > 0$ such that for every function 
$\psi\in\Mer(\faqd,\Hyp,d',V)$, we have 
$$
\sup_{\omega_0 + i \faqd} (1 + |\gl|)^{n} e^{-M|\!\Re(\gl)|}\|\psi(\gl)\|
\leq C 
\sup_{\omega + i \faqd} (1 + |\gl|)^{n}
e^{-M|\!\Re(\gl)|}\|q(\gl)p(\gl)\psi(\gl)\|. \Eq\Eqmb
$$
Let now $\gf \in \Mer(\faqd, \Hyp, d, V).$ Then
$\psi = \pi_{\omega_0, d} \gf  $ 
belongs to $\Mer(\faqd,\Hyp,d',V)$, 
so that (\Eqmb) holds.
This estimate remains valid if the supremum in the right hand side is
taken over the bigger set $\omega_1+i\faqd$. Since 
$\pi_{\omega_1, d} = q\pi_{\omega_0, d}$,
the required estimate (\Eqma) follows. 
\endProof

\Sect\SectC{The Fourier transform maps into $\PWd$}
We have already seen that $\Fou f\in
\Mer(\faqd,\Hyp(X,\tau),d_{X,\tau})\otimes\oC$
for $f\in\distc$. 
In order to show that $\varphi=\Fou f$ belongs to the Paley-Wiener space
we must verify both the estimate (\Eqbd) for 
some $n\in\cN(\faqd)$, 
and the conditions (i) and (ii) of Definition \DefBE{}.

Let $f\in\distM$.
The following estimate for $\Fou f$, from which both
(\Eqbd) and (ii) follow easily by application of
Lemma \LemID,
will now be established.
Let $R\in\real$, then there exists a polynomial $p\in\Pi_\Sigma(\faqd)$ 
and a constant $n\in\Natural$
such that
$$\sup_{\gl\in\baraqd(P,R)} (1+|\gl|)^{-n} e^{-M\,|\!\Re\gl|}
\|p(\gl)\Fou f(\gl)\|<\infty.\Eq\Eqcb$$

The verification of (\Eqcb) is based on the following
estimate for the Eisenstein integral (cf.\ [\rBSFI],
Lemma 4.3).
There exists a polynomial $p\in\Pi_\Sigma(\faqd)$ 
and for each $u\in U(\fg)$
a constant $n\in\Natural$ such that 
$$\sup_{x\in X_M, \gl\in\baraqd(P,R)} 
(1+|\gl|)^{-n} e^{-M\,|\!\Re\gl|} \|p(\gl)\East(\gl\col u;x)\| 
<\infty\Eq\Eqca$$
for all $M>0$.

Let $\tau_X\colon X\to\real$ be the map defined by
$\tau_X(kaH)=\|\log a\|$ for $k\in K$, $a\in A_q$.
It is easily seen that $\tau_X(k\exp Y\, H)=\|Y\|$ for
$k\in K$, $Y\in\fp\cap\fq$, hence it follows from
[\rSbook], Prop.\ 7.1.2,
that $\tau_X$ is smooth on the open subset of
$X$ where $\tau_X>0$. For $x\in X$ we have $x\in X_M$ if and only if
$\tau_X(x)\leq M$.

Let $h\in\Cinf(\real)$ be an arbitrary smooth function
satisfying $h(s)=1$ for $s\leq\frac12$ and $h(s)=0$ for $s\geq1$.
Then the function 
$$\varphi_\gl(x)=h(|\gl|(\tau_X(x)-M))
p(\gl)\East(\gl\col x)$$
is smooth and coincides with $p(\gl)\East(\gl\col x)$
in a neighborhood of $X_M$. Hence
$$p(\gl)\Fou f(\gl)= \int_X \varphi_\gl(x) f(x)\, dx,\Eq\Eqcc$$
and hence by (\Eqyb) 
$$\|p(\gl)\Fou f(\gl)\|\leq C \sup_{|\alpha|\leq k, x\in X}
\|L_{X^\alpha}\varphi_\lambda(x)\|$$
with constants $C$ and $k$ independent of $\gl$.
It follows from the Leibniz rule that
$$\|L_{X^\alpha}\varphi_\lambda(x)\|$$
is bounded by a constant times 
the product of
$$\sup_{|\beta|\leq k}  |h(|\gl|(\tau_X(X^\beta;x)-M))|$$
and
$$\sup_{|\beta|\leq k}\|p(\gl)\East(\gl\col X^\beta;x)\|.$$
The former factor is bounded by a constant times
$(1+|\gl|)^k$ and it vanishes outside
$X_{M+|\gl|^{-1}}$. By (\Eqca)
the second factor is estimated on this set by 
a constant times
$$(1+|\gl|)^{n} e^{(M+|\gl|^{-1})|\!\Re\gl|}\leq
(1+|\gl|)^{n} e^{1+M|\!\Re\gl|}$$
so that the
desired estimate (\Eqcb) follows.

It remains to be established that $\cL\Fou f=0$ for $\cL\in\ACR$.
Recall that  
$\East(\gl\col \,\cdot\,)$
is meromorphic in $\lambda$ 
with values in $C^\infty(X\col\tau^*)\otimes\oC$,
and that by definition 
an element $\cL\in\Mer(\faqcd,\Sigma)^*_\laur\otimes\oC^*$
with real support
belongs to $\ACR$ if and only if it annihilates $\gl\mapsto
\East(\gl\col\,\cdot\,)$.
The Fourier transform $\Fou f(\gl)\in\oC$ is obtained
by applying the linear form $f\,dx\in C^\infty(X\col\tau^*)'$ 
to
$\East(\gl\col \,\cdot\,)\in C^\infty(X\col\tau^*)\otimes\oC$.
We claim that the applications of $\cL$ and $f\,dx$
commute, so that
$$\cL\Fou f=\cL(f\,dx(\East(\,\cdot\,\col\,\cdot\,))
=f\,dx(\cL\East(\,\cdot\,\col\,\cdot\,))=0.$$
This claim is easily verified with the 
lemma below. 

\Lem\LemCA Let $\cL\in\Mer(\faqdc,\Sigma)^*_\laur$ be a
$\Sigma$-Laurent functional on $\faqdc$, and let
$\varphi\in \Mer(\faqdc,\Sigma,V)$, where $V$ is a complete
locally convex space. For each continuous linear form $\xi$
on $V$, the function $\xi\circ\varphi$ belongs to
$\Mer(\faqdc,\Sigma)$ and the following identity holds
$$\cL(\xi\circ\varphi)=\xi(\cL\varphi).$$
\endLem

\Proof We refer to [\rBSaf], Section 10, for notation.
We may assume that $\cL$ is supported in a single
point $a\in\faqdc$. If $\psi\in \cO_a(\faqdc,V)$ then
$\xi\circ\psi\in \cO_a(\faqdc)$ and $u(\xi\circ\psi)=
\xi(u\psi)$ for $u\in S(\faqdc)$.
The proof is now straightforward
from [\rBSaf], Definition 10.1 (see also Remark~10.2).
\endProof

\Sect\SectD{Distributional wave packets}
Recall that if $\gf\colon i\faqd\to\oC$ is continuous
and satisfies  the estimate
$$\sup_{\gl\in i\faqd} (1+|\gl|)^n\|\gf(\gl)\| <\infty\Eq\Eqdc$$
for each $n\in\Natural$, then we define the wave packet
$\Wave\gf\in\Cinf(X\col\tau)$ by
$$\Wave\gf(x)= \int_{i\faqd} \nE(\gf(\gl)\col\gl\col x)\,d\gl.
\Eq\Eqdd$$
The wave packet is related to the Fourier transform by
$$\inp{\Wave \gf}{g} = \inp{\gf}{\Fou g}
$$
for all $g\in\CciXt$, that is,
$$\int_X \inp{\Wave\gf(x)}{g(x)}\,dx=
\int_{i\faqd} \inp{\gf(\gl)}{\Fou g(\gl)}\,d\gl.\Eq\Eqda$$
The brackets in the latter equation refer to the 
sesqui-linear inner products on the 
finite dimensional Hilbert spaces 
$V_\tau$ and $\oC$, respectively.

The transform $\Wave$ can
be extended as follows to all continuous functions
$\gf\colon i\faqd\to\oC$ satisfying an estimate
$$\sup_{\gl\in i\faqd} (1+|\gl|)^{-n}\|\gf(\gl)\| <\infty,
\Eq\Eqdb$$
for some $n\in\Natural$. For such a function $\gf$
we define the distributional wave packet 
$\Wave\gf\in\Cmi(X\col\tau)$ by requiring (\Eqda) 
for all $g\in\CciXt$. It follows from the 
estimate (\Eqdb) together with  Lemma \LemBA,
that the integral on the right hand side of (\Eqda) is well-defined and
depends continuously on $g$, so that an element in
$\Cmi(X\col\tau)$ is defined by this equation. 

In particular, since for each $\gf\in\PWd$ the restriction
$\gf|_{i\faqd}$ is well-defined and satisfies
(\Eqdb) for some $n$, we thus have a well
defined linear map $$\Wave\colon\PWd\to\Cmi(X\col\tau).$$

\Sect\SectE{The Fourier transform is injective}

The injectivity is established in the following theorem.

\Thm\ThmEA There exists an invariant differential operator
$D\in\DGH$ which is formally selfadjoint,
injective as an operator 
$\Ccmi(X)\to\Ccmi(X)$
and which satisfies 
$$D\Wave\Fou f=\Wave\Fou Df=Df\Eq\Eqea$$
for all $f\in\distc$. 

In particular, the Fourier transform $\Fou\colon
\Ccmi(X\col\tau)\to \PWd$ is injective. 
\endThm

The proof will be given after the following lemma.

\Lem\LemEC Let $\varphi\colon i\faqd\to\oC$
be a continuous function satisfying {\rm (\Eqdc)} for all $n\in\Natural$.
Then  
$$
\int_X \inp{f(x)}{\Wave\gf(x)}\,dx
=
\int_{i\faqd} \inp{\Fou f(\gl)}{\gf (\gl)}\,d\gl.\Eq\Eqeb$$
for all $f\in\distc$.
\endLem

\Proof 
By taking adjoints in the estimate (\Eqca) a similar estimate
is derived for the Eisenstein integral $\nE(\gl\col x)$ and
its derivatives with respect to $x$.
It follows that 
the definition (\Eqdd) of $\Wave\gf$ allows an interpretation 
as an
integral over $i\faqd$
with values in the Fr\'echet space $\Cinf(X\col\tau)$.
In the left hand side of (\Eqeb) we apply $f\,dx$ to $\Wave\gf$.
By continuity we may then take $f\,dx$ 
inside the integral over $i\faqd$
and obtain
$$
\int_X \inp{f(x)}{\Wave\gf(x)}\,dx
=
\int_{i\faqd}\int_X 
\inp{f(x)}{\nE(\gf(\gl)\col\gl\col x)}\,dx\,d\gl,$$
which, by definition of $\East(\gl\col x)$,
exactly equals the right hand side of 
(\Eqeb).
\endProof

\Proofof{Theorem \ThmEA}
It follows from [\rBSmc] Thm.\ 14.1, Prop.\ 15.2 and 
Lemma 15.3
that there exists an invariant differential operator
$D\in\DGH$ which is formally selfadjoint,
injective as an operator $\Cci(X)\to
\Cci(X)$ and which satisfies 
(\Eqea)
for all $f\in\CciXt$. Since $D=D^*$, one obtains (\Eqea)
for $f\in\distc$ by transposition using (\Eqda) and (\Eqeb). 
Finally, it follows from Lemma \LemEB{} below,
that $D$ is injective on
the space of generalized functions as well. 
The injectivity of $\Fou$ is an immediate consequence
of (\Eqea) and the injectivity of $D$.
\endProof

\Lem\LemEB Let $D\in\DGH$. If $D$ is injective  $\Cci(X)\to
\Cci(X)$ then $D$ is injective $\Ccmi(X)\to\Ccmi(X)$.
\endLem

\Proof Recall that for 
$\phi\in \Cci(G)$ and $f\in\Cmi(X)$ we define
$L(\phi)f\in\Cinf(X)$ by
$$L(\phi)f(x)=\int_G \phi(g)f(g^{-1}x)\,dg.$$
The integral can be interpreted as a $\Cmi(X)$-valued
integral  
in the variable $x$, or it can be defined as the transpose
of the operator $L(\phi^\vee)\colon \Cci(X)\to\Cci(X)$.
In any case, $L(\phi)f$ is a smooth function on $G$
and it is compactly supported when $f$ has compact support.
Furthermore, $L(\phi)$ commutes with every invariant differential
operator $D$.

Let $\phi_j\in \Cci(G)$, $j\in\Natural$,  be an 
approximative unit, then it is well known that
$L(\phi_j)f$ converges weakly (in fact, also strongly)
to $f$, for each $f\in\Cmi(X)$. 

After these preparations the proof of the lemma is simple.
If $f\in\Ccmi(X)$ and $Df=0$ then 
$D(L(\phi_j)f)=L(\phi_j)Df=0$ and hence $L(\phi_j)f=0$ for all $j$. 
Hence 
$f=0$.
\endProof 

\Sect\SectL{Generalized Eisenstein integrals and 
Fourier transforms} 

Let $F\subset\Delta$, where $\Delta$ is the set of simple roots 
for $\Sigma^+$. We will use the notation of [\rBSPW], Section 5
and [\rBSPlI], Section 9.
In particular, $\CF$ is the finite dimensional Hilbert
space and
$\nE_F(\nu\col x)\in\Hom(\CF,V_\tau)$
the generalized Eisenstein integral defined in
eqs.\ (5.5)-(5.6) of [\rBSPW],
for $\nu\in\faqFdc$ and $x\in X$. 
The generalized Eisenstein integral  is a meromorphic 
$\Hom(\CF,V_\tau)$-valued
function of $\nu$, with singularities along a
real $\Sigma_r(F)$-configuration of
hyperplanes in $\faqFdc$ (see [\rBSPlI], Lemma 9.8). Here
$\Sigma_r(F)$ is the set of all non-zero restrictions
to $\faqF$ of elements in $\Sigma$.
For $F=\emptyset$ the generalized Eisenstein integral $\nE_F(\nu\col
x)$
is
identical with the normalized
Eisenstein integral $\nE(\gl\col x)$.

The corresponding {\it generalized Fourier transform} 
%$\Fou_F\colon \CciXt\to\Mer(\faqFdc,\Hyp_F)\otimes\Cusp_F$
is defined by
$$\Fou_F f(\nu)= \int_X \East_F(\nu\col x) f(x)\, dx\in\cA_F$$
for $\nu\in\faqFdc$, $f\in\CciXt$,
where $\East_F(\nu\col x)
=\nE_F(-\bar\nu\col x)^*\in\Hom(V_\tau,\Cusp_F)$.
The generalized Fourier transform is a meromorphic 
$\CF$-valued
function of $\nu$,

It is a remarkable property of
the generalized Eisenstein integral  $\nE_F(\nu\col x)$
that it can be obtained from the 
ordinary Eisenstein integral $\nE(\gl\col x)$ 
by applying a suitable operator in the
variable $\lambda$. More precisely, we have the following
result.

\Lem\LemGA There exists a Laurent functional
$\cL\in\Mer(\faqFdperpc,\gS_F)^*_\laur\otimes\Hom(\CF,\oC)$
with real support, such that
$$\nE_F(\nu\col x)=\cL[\nE(\nu+\,\cdot\,\col x)]$$
for $\nu\in\faqFdc$. 
\endLem

\Proof This follows immediately from [\rBSPlI] Lemma 9.7 with
$\psi$ in a basis for the finite dimensional
space $\CF$.\endProof

By taking adjoints it follows from Lemma \LemGA{}  
that there exists a 
Laurent functional 
$\cL_F\in\Mer(\faqFdperpc,\gS_F)^*_\laur\otimes\Hom(\oC,\CF)$
with real support such that
$$\East_F(\nu\col x)
=\cL_F[\East(\nu+\,\cdot\,\col x)].$$
Let such a Laurent functional, denoted by $\cL_F$, be fixed
in the sequel. It follows immediately that
$$\Fou_F f(\nu)=\cL_F[\Fou f(\nu+\,\cdot\,)].\Eq\Eqgb$$
In order to study the consequences of (\Eqgb) for
generalized Fourier transform, we need
the following result.

\Lem\LemHA Let $\Hyp$ be a real $\gS$-configuration in $\faqdc$,
and let $F\subset\Delta$ be given.
Let $\cL\in\Mer(\faqFdperpc,\gS_F)^*_\laur$ have
real support. There exists a real $\gS_r(F)$-configuration
$\Hyp_F$ in $\faqFdc$
and for every map $d\colon\Hyp \to \Natural$ a map $d'\colon
\Hyp_F \to \Natural$ such that the following holds
\smallskip\noindent{\rm (a)}
The operator $\cL_*$ defined by
$\cL_*\psi(\nu)= \cL[\psi(\nu+\cdot)]$ for
$\nu\in\faqFdc$ 
maps $\Mer(\faqd, \Hyp, d)$
continuously into $\Mer(\faqFd, \Hyp_F, d').$
\smallskip
\noindent{\rm (b)}
The operator $\Lau_*$ restricts to a continuous linear map
$
\Par(\faqd, \Hyp, d)\to \Par(\faqFd, \Hyp_F, d'),
$
where $\Par(\faqd, \Hyp, d)$ and
$\Pal(\faqFd,\Hyp_F,d')$ are defined in Definition \DefBA.
%\smallskip\noindent{\rm (c)}
%For every
%$n \in \Natural$ the operator $\Lau_*$ restricts to a continuous
%linear map
%$$
%\Par^*(\faqd,\Hyp,d,n) \longrightarrow
%\Par^*(\faqFd,\Hyp_F,d',n).
%$$
\endLem

\Proof Let $\Hyp_F$ and $d'$ be as in [\rBSaf], Cor.\ 11.6 (b).
Then $\cL_*$ maps $\Mer(\faqd,\Hyp,d)$ continuously
into $\Mer(\faqFd,\Hyp_F,d')$.
The assertion (b) is given in [\rBSPW], Lemma 6.1(v), with a proof
following  [\rBSrc] Lemmas 1.10-1.11.
%The proof of assertion (c) is similar.
\endProof

We fix a $\Sigma_r(F)$-hyperplane configuration 
$\Hyp(X,\tau,F)$ as $\Hyp_F$ in 
Lemma \LemHA, where we take $\Hyp=\Hyp(X,\tau)$ and $\Lau=\Lau_F$.
In addition, we fix a map $d_{X,\tau,F}\colon \Hyp(X,\tau,F) \to \Natural$
as $d'$, where we take $d=d_{X,\tau}$.
It follows from (a) that 
$$\East_F(\nu\col x)\in \Mer(\faqFd,\Hyp(X,\tau,F),
d_{X,\tau,F})\otimes\Hom(V_\tau,\CF).$$
Furthermore, the following result is obtained from (b), (\Eqgb)
and Lemma \LemBA.

\Lem\LemJA The generalized Fourier transform
is continuous
$$
\Fou_F\colon\CciXt\to \Pal(\faqFd,\Hyp(X,\tau,F),d_{X,\tau,F})\otimes\cA_F.$$
\endLem

\Sect\SectG{Generalized wave packets and
Fourier inversion}

Let $\Hyp_F=\Hyp(X,\tau,F)$. 
For
$\gf\in\Pal(\faqFd,\Hyp_F)\otimes\cA_F$
we introduce the {\it generalized wave packet}
$$\Wave_F\varphi(x)=\int_{\epsilon_F+i\faqFd}
\nE_F(\nu\col x)\gf(\nu)\,d\mu_{\faqdF}(\nu),
\Eq\Eqgd$$
where the element $\geps_F\in\fa^{*+}_{F\iq}$ and the measure
$d\mu_{\faqdF}$ on $\geps_F+i\faqFd$ 
are as defined in [\rBSFI], p.\ 42.
The definition is justified by the following estimate,
for which we refer to [\rBSFI], Lemma 10.8.
Let $\omega\subset \fa^{*}_{F\iq}$ be compact.  
There exists a polynomial on $\fa_{F\iq}^*$, $p\in\Pi(\fa_{F\iq})$,
for each $u\in U(\fg)$ a number $n\in\Natural$,
and for each $x$ a constant $C$,
locally uniform in $x$, such that 
$$\|p(\nu)\nE_F(\nu\col u; x)\|\leq C(1+|\nu|)^n
\Eq\Eqgc$$
for all $\nu\in\omega +i \faqFd$
(see also [\rBSPlI], Prop.\ 18.10, where a stronger result is given).
It follows  that for $\geps_F\in \fa^{*+}_{F\iq}$ 
sufficiently close to $0$, 
the integral (\Eqgd) is independent of $\geps_F$ and  
converges locally uniformly in $x$. Moreover, the resulting function
$\Wave_F\gf$ belongs to  $\Cinf(X\col\tau)$.

The Fourier inversion formula of [\rBSFI] now takes the form 
$$f(x)=\sum_{F\subset\Delta} c_F\int_{\epsilon_F+i\faqFd}
\nE_F(\nu\col x)\cL_F[\Fou f(\nu+\,\cdot\,)]
\,d\mu_{\faqdF}(\nu)=
\sum_{F\subset\Delta} c_F\Wave_F(\cL_{F*}\Fou f)(x)\Eq\Eqga
$$
for $f\in\CciXt$ (see [\rBSPW], Thm.\ 8.3), where the asterisk on
$\cL_F$ as in Lemma \LemHA{} indicates that it acts by
$$\cL_{F*}\varphi(\nu)=\cL_F[\varphi(\nu+\,\cdot\,)]$$
for $\gf\in\Mer(\faqd,\Hyp)$. 
The $c_F$ are explicitly given
constants. The explicit expression $c_F=|W| t(\faqF^+)$
is not relevant for the proof of Theorem \ThmBF.

\Sect\SectJ{Generalized distribution wave packets}
We shall see later that the inversion formula (\Eqga)
is valid also for distributions.
For this purpose we need to extend the
generalized wave packet map
$\Wave_F \colon
\Par(\faqFd,\Hyp_F)\otimes\cA_F\to\Cinf(X\col\tau)$ to a map
$\Par^*(\faqFd,\Hyp_F)\otimes \cA_F\to \Cmi(X\col\tau)$.
Here $\Par^*(\faqFd,\Hyp_F)$ is 
defined 
as $\Par^*(\faqd,\Hyp)$
in Definition \DefBD,
but with $\faqFd$ in place of $\faqd$  and 
$\Hyp_F=\Hyp(X,\tau,F)$
in place of~$\Hyp$.

In analogy with the definition 
of the distributional wave packet $\Wave\varphi$
given in 
Section~\SectD, we use an adjoint relation
with $\Fou_F$. For this purpose 
we introduce the sesquilinear pairing
$$
\Pal(\faqFd, \Hyp_F) \times \Pal^*(\faqFd, \Hyp_F) \to \complex
$$
given by
$$\inp{\varphi}{\psi}_\epsilon
=\int_{\epsilon_F+i\faqFd} \inp{\varphi(\gl)}{\psi(-\bar\gl)}\,
d\mu_{\faqdF}(\gl)\Eq\Eqjc$$  
for $\epsilon_F \in {\fa}_{Fq}^{*+}$ sufficiently close to zero.
The pairing $\inp{\,\cdot\,}{\,\cdot\,}$ inside the
integral is the standard sesquilinear pairing 
$\complex\times\complex\to\complex$.
The condition on $\epsilon_F$ guarantees that the domain of
integration is 
disjoint
from the singular locus of the integrand; moreover, by Cauchy's theorem
the integral is independent of the precise location of $\epsilon_F.$

For $d\colon \Hyp_F\to\Natural$ and $n\in\cN(\faqFd)$
the spaces $\Pal(\faqF, \Hyp_F, d)$ and 
$\Pal^*(\faqF, \Hyp_F, d, n)$ are defined and topologized as in
Section~\SectB.
The following lemma is obvious from these definitions.

\Lem\LemExtra
For each pair $d_1, d_2\colon\Hyp_F\to\Natural$ 
and every $n \in \cN(\faqFd),$ the pairing
{\rm (\Eqjc)}
restricts to a continuous sesquilinear pairing of 

%Fr\'echet spaces
$\Pal(\faqF, \Hyp_F, d_1)$ and $\Pal^*(\faqF, \Hyp_F, d_2, n).$
\endLem

The pairing (\Eqjc) is extended to
$\cA_F$-valued functions on $\faqFdc$ in the obvious fashion.
It is then easily seen by Fubini's theorem and (\Eqgc) that
$$\inp{\Wave_F\gf}{f}=\inp{\gf}{\Fou_F f}_\epsilon
\Eq\Eqja$$
for $\gf\in\Par(\faqFd,\Hyp_F)\otimes\cA_F$, $f\in\CciXt$.

The {\it generalized distribution wave packet} $\Wave_F\gf\in
\Cmi(X\col\tau)$ is defined for functions $\gf\in 
\Par^*(\faqFd,\Hyp_F)\otimes\cA_F$
by (\Eqja) for all $f\in\CciXt$.
The definition is justified by Lemmas \LemJA{} and
\LemExtra.

The following lemma is immediate from the definition.
 
\Lem\LemJB Let $\Hyp_F=\Hyp(X,\tau,F)$.
Let $d\colon\Hyp_F\to\Natural$ and $n \in \cN(\faqFd),$ 
be arbitrary.
The distributional generalized wave packet map
$$\Wave_F\colon\Pal^*(\faqFd, \Hyp_F, d,n)\otimes\cA_F 
\to C^{-\infty}(X\col \tau)$$
is continuous for the weak topology with respect to the pairing 
{\rm (\Eqjc)} on 
$\Pal^*(\faqFd, \Hyp_F, d,n)$,
and the weak dual topology on $C^{-\infty}(X\col \tau)$.
%(b) the original topology on the first space
%and the strong dual topology on the second space.
\endLem

It follows from Lemma \LemExtra{} that the map is 
continuous also for the original topology on 
$\Pal^*(\faqFd, \Hyp_F, d,n)$ and the weak dual topology on 
$C^{-\infty}(X\col \tau)$. We shall see later
(in Lemma \LemIA) that it is continuous for the original  
topology on 
$\Pal^*(\faqFd, \Hyp_F, d,n)$ and the strong dual topology on 
$C^{-\infty}(X\col \tau)$. However,
this is not needed for the proof of Theorem~\ThmBF.

\Sect\SectF{Multiplication operators on the Paley-Wiener space}
The main result of this subsection, Prop.\ \PropFA, was 
announced jointly with Flensted-Jensen in 
the survey paper [\rBFJS], see Prop.\ 18.
However, the present proof is independent of 
the preceding results of that paper. 

Let $\fb$ be a Cartan subspace of $\fq$ containing
$\faq$. Then $\fb=\fbk\oplus\faq$
with $\fbk=\fb\cap\fk$. Let $\Wb$ be the Weyl group
of the restricted root system of $\fbc$ in $\fgc$.
Let $\fb^d$ denote the real form $i\fbk\oplus\faq$
of $\fbc$.

Let $\cO(\fbc^*)^\Wb$ denote the space of  $\Wb$-invariant entire 
functions on $\fbc^*$, and for $r>0$ let 
$\PW_r(\fb^d)^\Wb$ denote the subspace of functions which are also
rapidly decreasing of exponential type $r$,
that is, 
$$\sup_{\gl\in\fbc^*} (1+|\gl|)^n e^{-r|\!\Re\gl|}|\psi(\gl)|<\infty$$
for all $n\in\Natural$. The real part $\Re\gl$ is taken with respect
to the decomposition $\fbc=\fb^d+i\fb^d$.
Let $\PW(\fb^d)^\Wb$ denote the
union over $r$ of all the spaces $\PW_r(\fb^d)^\Wb$.

Given $\psi\in\cO(\fbc^*)^\Wb$ we define a multiplication operator
$M(\psi)$ on $\Mer(\faqdc,\oC)$ as follows. Recall
the orthogonal decomposition 
$$\oC=\oplus_{\gL\in L} \oC[\gL]\Eq\Eqfa$$
(see [\rBSmc], eq.\ (5.14)) where $L\subset i\fbk^*$ is a finite
set depending on $\tau$.
We define, for each $\gl\in\faqdc$, an endomorphism $M(\psi,\gl)$
of $\oC$ 
by $M(\psi,\gl)\eta=\psi(\gl+\gL)\eta$ for $\eta\in\oC[\gL]$,
and we define for each $\varphi\in\Mer(\faqdc,\oC)$ 
a function $M(\psi)\varphi\in\Mer(\faqdc,\oC)$ by 
$$M(\psi)\varphi(\gl)=M(\psi,\gl)\varphi(\gl)$$
for all $\gl\in\faqdc$.

The motivation behind this definition is as follows. 
If $\psi$ belongs to $\PW_r(\fb^d)^\Wb$ then the operator
$M(\psi)$ on $\Mer(\faqdc,\oC)$ corresponds, via the
Fourier transform $\Fou$, to a linear operator $M_\psi$ on 
$\Cci(X\col\tau)$, a so-called {\it multiplier}, so that
$$\Fou(M_\psi f)=M(\psi)\Fou f$$
for all $f\in\Cci(X\col\tau)$. The existence of the multiplier
$M_\psi$ is given a relatively elementary proof in [\rBFJS] 
without reference to
the Paley-Wiener theorem for $\Cci(X\col\tau)$, which was only
a conjecture when that paper was written. However, with the
Paley-Wiener theorem for $\Cci(X\col\tau)$ available from [\rBSPW],
the existence of the multiplier is an immediate consequence
of Prop.\ \PropFA{} below. 

Let $D\in\DGH$.
It follows from [\rBSmc], Lemma 6.2, that 
$$\Fou(Df)(\gl)=\mu(D,\gl)\Fou f(\gl).$$
By the definition of the decomposition
(\Eqfa), the endomorphism $\mu(D,\gl)$ of $\oC$ acts on
$\oC[\gL]$ as multiplication with $\gamma(D,\gl+\gL)$
where 
$$\gamma\colon\DGH\to P(\fb^*_\iC)^\Wb$$ is the
Harish-Chandra isomorphism.
Hence
$$\Fou(Df)=M(\gamma(D))\Fou f\Eq\Eqfb$$
for $f\in\CciXt$.

\Prop\PropFA Let $\psi\in\PW(\fb^d)^\Wb$.
The multiplication operator $M(\psi)$ maps the space
$\PWd$ into $\PWXt$.
More precisely, if $r,R>0$ and $\psi\in\PW_r(\fb^d)^\Wb$, then
$M(\psi)$ maps
$\PW^*_R(X\col\tau)$ into $\PW_{R+r}(X\col\tau)$.
\endProp

\Proof The idea of the proof is taken from [\rArth], p.\ 87.
Let $\psi\in\PW_r(\fb^d)^\Wb$, then since $L$ is finite 
there exists, 
for each $n\in\Natural$ a constant
$C>0$ such that 
$$
|\psi(\lambda+\Lambda)|\leq C(1+|\gl|)^{-n}e^{r|\Re\gl|}
$$
for all $\gl\in\faqcd$ and $\Lambda\in L$.
It is now easily seen that the estimates 
(\Eqba) and (ii)
in 
Definition \DefBC{}
of $\PW_{R+r}(X\col\tau)$ are satisfied
by $M(\psi)\gf$ for $\gf\in \PW_R^*(X\col\tau)$. Only
the annihilation by $\ACR$ in
Definition \DefBC{} (i) remains to be verified.

Let $\gf\in\Mer(\faqdc,\gS,\oC)$ 
be a function annihilated by 
all $\cL\in\ACR$, and let $\psi\in\cO(\fbc^*)^\Wb$.
We claim that then $M(\psi)\gf$
is also annihilated by all  $\cL\in\ACR$.

Let  $\cL\in\ACR$ and
assume first that $\psi$ is a polynomial. Then there exists
an invariant differential operator $D\in\DGH$ such that
$\psi=\gamma(D)$. It follows from (\Eqfb) that 
$\cL(M(\psi)\Fou f)=\cL(\Fou(Df))=0$ for all $f\in\CciXt$.
By [\rBSaf], p.\ 674, there exists a Laurent
functional $\cL'\in\LauoC$
such that $\cL(M(\psi)\phi)=\cL'\phi$
for all $\phi\in\Mer(\faqdc,\gS,\oC)$. 
Moreover, $\supp\cL'\subset\supp\cL$.
Hence $\cL'\in\ACR$ 
by [\rBSPW], Lemma 3.8,
and we conclude that
$\cL(M(\psi)\gf)=\cL'\gf=0$.

Consider now the case of a general function $\psi\in\cO(\fbc^*)^\Wb$.
We expand $\psi$ in its Taylor series around $0$, and denote by
$\psi_k$ the sum of the terms up to degree $k$. Then $\psi_k\to\psi$,
uniformly on compact sets, from which it follows that
$\cL(M(\psi_k)\gf)\to\cL(M(\psi)\gf)$.
Each $\psi_k$ is a $\Wb$-invariant polynomial,
hence $\cL(M(\psi_k)\gf)=0$. It follows that 
$\cL(M(\psi)\gf)=0$.
\endProof

\Sect\SectH{The Fourier transform is surjective}
Let $\varphi\in\PWMd$ be given. Inspired by Equation (\Eqga)
we define
$$f=\sum_{F\subset\Delta}
c_F\Wave_F(\cL_{F*}\varphi)\in\Cmi(X\col\tau),\Eq\Eqhb$$
where $\cL_F$ is chosen as in Section \SectL, see (\Eqgb). 
For (\Eqhb) to make sense we need that
$\cL_{F*}\varphi$ belongs to the space
$\Par^*(\faqFd,\Hyp_F)\otimes\cA_F$ on which $\Wave_F$
was defined, see (\Eqja). This is secured by the following 
lemma.

\Lem\LemHB Let $\Hyp$, 
$\cL$,  
$\Hyp_F$, $d$ and $d'$ be as in Lemma \LemHA.
For every $n \in \cN(\faqd)$ there exists $n'\in \cN(\faqFd)$
such that
the operator $\Lau_*$ restricts to a continuous
linear map
$$
\Par^*(\faqd,\Hyp,d,n) \longrightarrow
\Par^*(\faqFd,\Hyp_F,d',n').
$$
\endLem

\Proof 
The proof is similar to the proof of Lemma \LemHA{} (b).
\endProof

It follows that (\Eqhb) defines a distribution 
$f\in\Cmi(X\col\tau)$.
We claim that $f\in\Cmi_M(X\col\tau)$ 
and $\Fou f=\varphi.$

Let $\psi_j\in\PW(\fb^d)^W_{r_j}$ be a sequence of functions
such that $r_j\to 0$ for $j\to\infty$, such that
$\psi_j$ is uniformly bounded on each set of the form
$\omega+i\fb^{d*}$ with $\omega\subset\fb^{d*}$ compact,
and such that $\psi_j\to 1$, locally uniformly on $\faqdc.$
Such a  sequence
can be constructed by application of the Euclidean Fourier
transform to a smooth approximation of the Dirac measure on $\fb^d$.
In particular,  for each $\gL\in L$, the sequence of functions
$\psi_j(\,\cdot\,+\gL)$ is uniformly bounded on each set
$\omega+i\fb^{d*}$ as above, and converges to $1,$ locally uniformly
on~$\faqdc.$

Consider the functions $\varphi_j: = M(\psi_j)\varphi$. It follows from
Proposition \PropFA{} that $\varphi_j\in
\PW_{M+r_j}(X\col\tau).$ Hence, by the Paley-Wiener
theorem of [\rBSPW] there exists a unique function $f_j
\in \Cci(X\col\tau)$, with support in $X_{M+r_j}$,
such that $\Fou f_j=M(\psi_j)\varphi$.

Let $d$, $n$ be such that $\varphi\in \Par^*(\faqd,\Hyp,d,n)
\otimes \oC$. Then it follows frow the properties of $\psi_j$ mentioned
above that $\cF f_j = \gf_j \to \gf$ for $j \to \infty$ as a sequence of
functions in $\Mer(\faqd, \Hyp, d).$
Moreover, the sequence is bounded as a sequence in $\Par^*(\faqd,\Hyp,d,n).$

In view of Lemmas \LemHA{} and \LemHB{}
it follows that the sequence $\cL_{F*}\Fou f_j$
is bounded in 
$\Pal^*(\faqFd,\Hyp_F,d',n')\otimes \cA_F$ 
%for some $n'\in\Natural$,
and that $\cL_{F*}\Fou f_j \to \cL_{F*}\varphi$ as a sequence
in $\Mer(\faqFd, \Hyp_F, d').$ By dominated convergence it follows
that $\cL_{F*}\Fou f_j \to \cL_{F*}\varphi$ weakly in 
$\Pal^*(\faqFd, \Hyp_F, d',n')$  
with respect to the pairing (\Eqjc).
In view of Lemma \LemJB{} this implies
that
$$
f_j=\sum_Fc_F\Wave_F\Lau_{F*}\Fou f_j\to
\sum_Fc_F\Wave_F\Lau_F\varphi=f
$$
weakly in $\Cmi(X\col\tau)$. Since $f_j$ belongs to
$\Cinf_{M+r_j}(X\col\tau)$
for each $j$, we conclude that $f\in\Cmi_M(X\col\tau)$.
Moreover, it follows from the weak convergence
$f_j\to f$ that $\Fou f_j(\gl)\to\Fou f(\gl)$
for all $\gl$ outside the hyperplanes in $\Hyp.$ Hence,
$\Fou f=\varphi$ as claimed. 

\smallskip
Theorem \ThmBF{} has now been proved.

\Cor\CorHA The Fourier inversion formula {\rm (\Eqga)}
$$f=\sum_{F\subset\Delta} c_F\Wave_F\cL_{F*}\Fou f\Eq\Eqha
$$
is valid for $f\in\Ccmi(X\col\tau)$.
\endCor

\Proof It was seen during the proof above that
$\sum_Fc_F\Wave_F\Lau_{F*}\varphi\in\Ccmi(X\col\tau)$ and
$$\varphi=\Fou(\sum_Fc_F\Wave_F\Lau_{F*}\varphi)$$
for all $\varphi\in\PWd$. In particular, the latter identity applies
to $\varphi=\Fou f$ for each $f\in\Ccmi(X\col\tau)$.
The formula (\Eqha) then follows from the injectivity of
$\Fou$ (Theorem \ThmEA).\endProof

\Sect\SectI{A topological Paley-Wiener theorem}
We shall equip the spaces $\distM$ and $\PWMd$ with
natural topologies for which the Fourier transform is an 
isomorphism.

On the space of generalized functions on $X$ 
we use the {\it strong dual topology}, where
we regard $C^{-\infty}(X)$ as the dual space of $\Cci(X)$.
Recall that by definition, the strong dual topology 
on $\Cci(X)^*$ is the locally convex topology 
given by the seminorm system
$$p_B(f)=\sup_{\varphi\in B} |f(\varphi)|$$
where $B$ belongs to the family of
all bounded subsets of $\Cci(X)$.
Notice that $\Cci(X)$ is a Montel space, that is, 
it is reflexive and a subset is bounded 
if and only if it is relatively compact
(see [\rSchaefer], p.\ 147). 

On the space $C^{-\infty}_c(X)$ of compactly supported
generalized functions on $X$ we use the strong dual topology, where
we regard $C^{-\infty}_c(X)$ as the dual space of $C^\infty(X)$.
As an immediate consequence of these dualities,
the inclusion map $C_c^{-\infty}(X)\to C^{-\infty}(X)$ is 
continuous, and multiplication by a function $\psi\in C_c^{\infty}$
is continuous $C^{-\infty}(X)\to C_c^{-\infty}(X)$.

The topologies on $C^{-\infty}(X)$ and $C^{-\infty}_c(X)$
induce the same topology on the space of distributions
supported in a fixed compact subset $\Omega$
of $X$. This follows from
the last remark of the preceding paragraph, when we take as $\psi$
any function which is identically 1 on a neighborhood of $\Omega$.

In particular, for 
each $M>0$ the space $\distM$ of $\tau$-spherical generalized 
functions with support in $X_M$ is topologized in this fashion,
as a topological subspace of $C^{-\infty}(X)\otimes V_\tau$,
or equivalently, as a topological subspace of 
$C_c^{-\infty}(X)\otimes V_\tau$.

Recall that $\PWMd\subset \Pal^*(\faqd,\Hyp(X,\tau),d_{X,\tau})
\otimes\oC$,
and that 
$$\Pal^*(\faqd,\Hyp,d)=\cup_{n\in\cN(\faqd)}\Pal^*(\faqd,\Hyp,d,n).$$
On each space $\Pal^*(\faqd,\Hyp,d,n)$, where
$d\colon\Hyp\to \Natural$ and $n\in\cN(\faqd)$, the topology was 
defined by means of the seminorms (\Eqbd). 
On $\cN(\faqd)$ we define an order relation by
$n_1\leq n_2$ if and only if $n_1(\omega)\leq n_2(\omega)$
for all $\omega$. It is easily seen that if $n_1\leq n_2$
then $\Pal^*(\faqd,\Hyp,d,n_1)\subset\Pal^*(\faqd,\Hyp,d,n_2)$
with continuous inclusion. The family of spaces  
$\Pal^*(\faqd,\Hyp,d,n)$ indexed by $n\in\cN$ is thus a directed family,
and we can give $\Pal^*(\faqd,\Hyp,d)$
the inductive limit topology for the union over $n$.
The Paley-Wiener space $\PWMd$ is given the relative
topology of this space (where $d=d_{X,\tau}$), tensored by
$\oC$.

\Thm\ThmKA The Fourier transform is a topological isomorphism of
$\distM$ onto $\PWMd$, for each $M>0$.\endThm

\Proof Only the topological statement remains to be proved.
We will prove that $\Fou$ is continuous with respect to the
topology of 
$C_c^{-\infty}(X)\otimes V_\tau$, and that its inverse
is continuous into $C^{-\infty}(X)\otimes V_\tau$.
Since the topologies agree on $\distM$, as remarked above, this will
prove the theorem.

Let $\Hyp=\Hyp(X,\tau)$ and $d=d_{X,\tau}$.
For the continuity of the Fourier transform
$$\Fou\colon C_c^{-\infty}(X\col\tau)
\to\Pal^*(\faqd,\Hyp,d) \otimes\oC
\Eq\Eqia$$ 
we 
remark that by a theorem
of Grothendieck, $C_c^{-\infty}$ is bornological,
since it is the strong dual of the reflexive Fr\'echet space
$C^{\infty}$ (see [\rSchaefer], p.\ 154). Therefore,
it suffices to prove that the Fourier transform maps every
bounded set $B\subset  C_c^{-\infty}(X\col\tau)$
to a bounded set in $\Pal^*(\faqd,\Hyp,d) \otimes\oC$
(see [\rSchaefer], p.\ 62). Since all such sets $B$ are 
equicontinuous (see [\rSchaefer], p.\ 127) we may assume 
that there exists a continuous seminorm $\nu$ on 
$C^\infty(X\col\tau^*)$ such that
$$B\subset\Big\{f\, \Big|\,\,\,\,  
\big|\!\int \varphi(x)f(x)\,dx\big|\leq \nu(\varphi),
\forall\varphi\in C^\infty(X\col\tau^*)\Big\}.\Eq\Eqkb$$

As in (\Eqyb) we may assume that the seminorm $\nu$ has the form
$$\nu(\varphi)=C\sup_{|\alpha|\leq k, x\in \Omega}
\| L_{X^\alpha}\varphi(x)\|$$
for some $C>0$, $k\in\Natural$ and $\Omega\subset X$ compact. 
Choose $M>0$ such that $\Omega\subset X_M$, and let 
$R\in\real$.
It follows from
(\Eqca) that there exists a number
$n\in\Natural$ 
and a polynomial $p\in\Pi_\gS(\faqd)$
such that
$$\sup_{\gl\in\baraqd(P,R)}(1+|\gl|)^{-n} e^{-M|\!\Re\gl|}
\nu(p(\gl)\East(\gl\col \,\cdot\,))<\infty.\Eq\Eqka$$
It now follows from (\Eqkb) with $\varphi=
p(\gl)\East(\gl\col\,\cdot\,)$,
combined with (\Eqka),
that 
$$\sup_{f\in B}\,\sup_{\gl\in\baraqd(P,R)} (1+|\gl|)^{-n} e^{-M\,|\!\Re\gl|}
\|p(\gl)\Fou f(\gl)\|<\infty.\Eq\Eqkc$$

For each compact set $\omega\subset\faqd$ we choose
$R\in\real$ such that $\omega\subset\baraqd(P,R)$ and define
$n(\omega)$ to be the number $n$ in (\Eqkc). 
By application of Lemma \LemID{} it follows that
the seminorm (\Eqbd) is uniformly bounded on $\Fou(B)$.
With $n\in\cN(\faqd)$ chosen in this fashion, we thus see that
$\Fou(B)$ is contained and
bounded in
$\Pal^*(\faqd,\Hyp,d,n) \otimes\oC$, hence also
in the union $\Pal^*(\faqd,\Hyp,d) \otimes\oC$
with the inductive limit topology. 
Thus (\Eqia) is continuous.

In order to establish the continuity of the inverse Fourier
transform we use Corollary \CorHA, 
according to which the inverse Fourier transform
is given by the finite sum  of $c_F$ times
$\Wave_F\cL_{F*}$.
The operator
$\cL_{F*}$ is continuous from 
$\Pal^*(\faqd,\Hyp(X,\tau),d_{X,\tau})\otimes \oC$ to
$\Pal^*(\faqFd,\Hyp(X,\tau,F),d_{X,\tau,F})\otimes \cA_F$ 
by 
Lemma \LemHB, and continuity of $\Wave_F$ 
is established in the lemma below.
\endProof

\Lem\LemIA The generalized wave packet operator $\Wave_F$
is strongly continuous 
$$\Pal^*(\faqFd,\Hyp(X,\tau,F),d_{X,\tau,F})\otimes \cA_F\to 
C^{-\infty}(X\col\tau)$$
for each $F\subset\Delta$. 
\endLem

\Proof
Let $B\subset \CciXt$ be bounded,
then $\Fou_F(B)$ is bounded by  Lemma \LemJA,
and it follows from  (\Eqja) that
$$p_B(\Wave_F\varphi)=\sup_{f\in B} 
\left|\inp{f}{\Wave_F\varphi}\right|
=\sup_{f\in B}
\left|\inp{\Fou_F f}{\varphi}_\epsilon\right|.
$$
Hence $p_B\circ\Wave_F$ is continuous  by Lemma \LemExtra. 
\endProof

\Sect\SectK{Further properties of the topology}

In this final section we show that $\Pal^*(\faqd,\Hyp,d,n)$
is a Fr\'echet space, under a certain natural condition on
$n\in\cN(\faqd)$, and that this condition is satisfied by
sufficiently many elements $n$ to give the same union
$\Pal^*(\faqd,\Hyp,d)$.

Let $\Hyp$ be a real $\gS$-configuration in $\faqdc$, and
let $d\colon\Hyp\to\Natural$ be arbitrary.
For $\omega\in\cC(\faqd)$ and $n\in\cN(\faqd)$, 
we denote by $\nu^*_{\omega,n}$ the seminorm on $\Pal^*(\faqd,\Hyp,d,n)$
defined in (\Eqbd).

\Def\DefIA A function $n\in\cN(\faqd)$ is called 
{\it regular} 
if for every $\omega\in\cC(\faqd),$ 
$$
n(\omega) = \inf \{ n(\omega') \mid \omega'\in\cC(\faqd), \,
\omega \subset {\rm int}(\omega') \}.
$$
The set of such functions $n$ is denoted $\cN_0(\faqd)$. 
\endDef

\smallskip
Notice that if $n$ is regular, then
$n(\omega_1)\leq n(\omega_2)$ for $\omega_1\subset\omega_2$.
It is also easy to see that if $n_1,n_2$ are regular, then so is
$n=\max(n_1,n_2)$.

\Lem\LemIB 
Let $n\in\cN_0(\faqd)$. Then $\cP^*(\faqd,\Hyp,d,n)$ is a Fr\' echet space.
\endLem

\Proof It is easily seen that $\cP_n^*$ 
is complete (also without the condition of regularity).
We have to show that it is metrizable. We will do this by
pointing out a countable collection of compact sets $\omega$,
such that the corresponding family of seminorms $\nu^*_{\omega,n}$
generates the topology.

The topology of $\faqd$ is locally compact and second countable.
Let $\cB$ be a countable basis of open sets with compact closures.
For every finite collection $F\subset\cB$, let
$$
\omega_F = \cup_{B \in F} \bar B.
$$ 
Then $\omega_F\in\cC(\faqd)$. We will show 
that the countable family
$
\nu^*_{\omega_F,n}
$ 
generates the topology.

Indeed, let $\omega\in\cC(\faqd)$. Then by regularity
there exists $\omega'\in\cC(\faqd)$,
with $n(\omega') = n(\omega),$ and $\omega\subset {\rm
int}(\omega')$. 
If $\omega''\in\cC(\faqd)$ satisfies $\omega\subset
\omega''\subset \omega'$,
it follows  that
$
n(\omega)=n(\omega'')=n(\omega').
$ 
We now see that we may take $\omega'$ so that in addition
$
\pi_{\omega', d} = \pi_{\omega, d}.
$

For every $\gl \in \omega$ there exists
$B\in\cB$ with $\gl \in B$ and $\bar B\subset {\rm int} 
(\omega').$ By compactness, there exists a finite collection
$F\subset\cB$ such 
that 
$
\omega\subset \omega_F \subset {\rm int }(\omega').
$
It follows that
$
n(\omega) =  n(\omega_F)
$
and $\pi_{\omega, d}=\pi_{\omega_F, d}.$
Hence 
$$\eqalign
{\nu^*_{\omega, n}(\gf) & = 
\sup_{\omega+ i \faqd} (1 + |\gl|)^{- n(\omega_F)} \|\pi_{\omega_F,
d}(\gl) \gf(\gl)\|\cr
& \leq 
\sup_{\omega_F+ i \faqd} (1 + |\gl|)^{- n(\omega_F)} \|\pi_{\omega_F,
d}(\gl) \gf(\gl)\| = \nu^*_{\omega_F,n}(\gf)}
$$
for all $\varphi$.\endProof

\Prop\PropIC For each $n\in\cN(\faqd)$ there exists a function
$\bar n\in\cN_0(\faqd)$ such that 
$$\cP^*(\faqd,\Hyp,d,n)\subset \cP^*(\faqd,\Hyp,d,\bar n)$$
with continuous inclusion.
\endProp

\Proof We define the map
$\bar n: \cC(\faqd) \to \Natural$ by
$$
\bar n(\omega) = \inf \{n(\omega')\mid \omega'\in\cC(\faqd), \,
\omega \subset {\rm int}(\omega')\}.
$$
It is easily seen that $\bar n$ is regular.

Let $\omega$ be a compact set. Then there exists a compact neighborhood
$\omega_1$ of $\omega$ such that $\bar n (\omega) = n(\omega_1).$
It follows from Lemma \LemID{} 
that there exists a constant $C > 0$ such that for
all $\gf \in \Mer(\faqdc, \Hyp, d, \oC),$
$$\nu^*_{\omega,\bar n}(\gf)=
\nu_{\omega, -\bar n(\omega)}(\gf) \leq C \nu_{\omega_1,-\bar n(\omega)}(\gf)
=
C\nu_{\omega_1,- n(\omega_1)}(\gf)=C\nu^*_{\omega_1, n}(\gf).
$$
The result follows.
\endProof

In the above we have proved that there exists a subset $\cN_0$
of $\cN = \cN(\faqd)$ with the following properties
(where $\cP^*_n=\cP^*(\faqd,\Hyp,d,n)$ for simplicity):

\smallskip
\item{(a)}
for every $n \in \cN_0$ the space $\cP^*_n$ is Fr\'echet;
\item{(b)}
for every $n_1, n_2 \in \cN_0$ there exists  $n \in \cN_0$
such that $n_1 \leq n,$ $n_2 \leq n$ (directed family);
\item{(c)}
for every $n \in \cN$ there exists $m \in \cN_0$
such that $\cP_n^* \subset \cP_m^*$ with continuous inclusion.
\smallskip

Because of these properties, the union $\cup_{n\in\cN_0} \cP^*_n$
is equal to $\cup_{n\in\cN} \cP^*_n=\cP^*(\faqd,\Hyp,d)$ and
the limit topologies defined by 
$\lim_{\cN_0} \cP^*_n$ and $\lim_{\cN} \cP^*_n$ 
are equal.
In particular, we see that $\cP^*(\faqd,\Hyp,d)$ is an inductive 
limit of Fr\'echet spaces (notice however, that it is not necessarily a
{\it strict} inductive limit).

\bigskip\goodbreak
\noindent{\bf References}
\def\Lable#1{\par\smallskip\noindent [#1]}

\Lable{\rArth} 
{\it J.\  Arthur,}
A Paley--Wiener theorem for real reductive groups,
Acta Math.\  {\bf 150} (1983), 1--89.

%\Lable{\rBJFA}
%{\it E.\ P.\ van den Ban,}
%The principal series for a reductive symmetric space, II. 
%Eisenstein integrals. 
%J.\ Funct.\ Anal.\ {\bf 109} (1992), 331-441.

\Lable{\rBsurvey}
{\it E.\ P.\ van den Ban,}
The Plancherel theorem for a reductive symmetric space, 
pp.\ 1-97 in
{\it Lie theory, Harmonic analysis on symmetric spaces -
General Plancherel theorems}, 
J.-P.\ Anker and B.\ \O rsted (eds.), Birkh\"auser Boston, 2005.

\Lable{\rBFJS}
{\it  E.\ P.\ van den Ban, M.\ Flensted-Jensen and H.\ Schlichtkrull,}
Basic harmonic analysis on pseudo-Riemannian symmetric
spaces, pp.\ 69-101 in {\it
Noncompact Lie groups and some of their applications},
E.\ A.\ Tanner and R.\ Wilson (eds.) (proc.\ of
the NATO advanced workshop, San Antonio, Texas 1993), Kluwer 1994.

%\Lable{\rBSfo} {\it E.\ P.\ van den Ban and H.\ Schlichtkrull,}
%Fourier transforms on a semisimple symmetric space.
%Invent.\ Math.\ {\bf 130} (1997), 517-574.

\Lable{\rBSmc} {\it E.\ P.\ van den Ban and H.\ Schlichtkrull,}
The most continuous part of the Plancherel decomposition
for a reductive symmetric space. 
Ann.\ Math.\ {\bf 145} (1997), 267-364.

\Lable{\rBSrc} {\it E.\ P.\ van den Ban and H.\ Schlichtkrull,}
A residue calculus for root systems. 
Compositio Math. {\bf 123} (2000), 27-72.

\Lable{\rBSFI} {\it E.\ P.\ van den Ban and H.\ Schlichtkrull,}
Fourier inversion on a reductive symmetric space.
Acta Math. {\bf 182} (1999), 25-85.

\Lable{\rBSaf} {\it E.\ P.\ van den Ban and H.\ Schlichtkrull,}
Analytic families of eigenfunctions on a
reductive symmetric space. Represent. Theory {\bf 5} (2001), 615-712.

\Lable{\rBSPlI} {\it E.\ P.\ van den Ban and H.\ Schlichtkrull,}
The Plancherel decomposition for a reductive symmetric space I.
Spherical functions. Invent.\ Math.\ {\bf 161} (2005), 453-566.

\Lable{\rBSPlII} {\it E.\ P.\ van den Ban and H.\ Schlichtkrull,}
The Plancherel decomposition for a reductive symmetric space II.
Representation theory.  Invent.\ Math.\ {\bf 161} (2005), 567-628.

\Lable{\rBSPW} {\it E.\ P.\ van den Ban and H.\ Schlichtkrull,}
A Paley-Wiener theorem for reductive symmetric
spaces. Preprint 2003, arXiv.math.RT/0302232.
%To appear in Ann.\ Math.

\Lable{\rBSvD} {\it E.\ P.\ van den Ban and H.\ Schlichtkrull,}
Paley-Wiener spaces for real reductive Lie groups. 
To appear in Indag.\ Math.

%\Lable{\rCamp} {\it O.\ A.\ Campoli}, 
%Paley-Wiener type theorems for rank-1 semisimple
%Lie groups,  Rev.\ Union Mat.\ Argent.\ {\bf 29} (1980), 197-221.

%\Lable{\rCD} {\it J.\ Carmona and P.\ Delorme,}
%Transformation de Fourier sur l'espace de Schwartz 
%d'un espace sym\'etrique r\'eductif.
%Invent.\ Math.\ {\bf 134} (1998), 59-99.

\Lable{\rDadok}
{\it J.\ Dadok,} Paley-Wiener theorem for singular support
of $K$-finite distributions on symmetric spaces.
J.\ Funct.\ Anal.\ {\bf 31} (1979), 341-354.

%\Lable{\rDel} {\it P.\ Delorme,}
%Formule de Plancherel pour les espaces sym\'etriques reductifs.
%Ann.\ of Math.\ {\bf 147} (1998), 417-452.

\Lable{\rEHO}
{\it M.\ Eguchi, M.\ Hashizume and K.\ Okamoto,}
The Paley-Wiener theorem for distributions
on symmetric spaces.
Hiroshima Math.\ J.\ {\bf 3} (1973), 109-120.

\Lable{\rGang} {\it R.\ Gangolli}, 
On the Plancherel formula and the Paley-Wiener theorem
for spherical functions on semisimple Lie groups, Ann.\ of Math.\ 
{\bf 93} (1971), 150-165.

%\Lable{\rHCIII} {\it Harish-Chandra,} 
%Harmonic analysis on real reductive
%groups III.\ The Maass-Selberg relations and the Plancherel
%formula, Ann.\ of Math.\ {\bf 104} (1976), 117-201.

\Lable{\rHePW} {\it  S.\ Helgason,}
An analogue of the Paley-Wiener theorem for the Fourier transform
on certain symmetric spaces, Math.\ Ann.\ {\bf 165} (1966), 297-308.

\Lable{\rHeBAMS} {\it  S.\ Helgason,}
Paley-Wiener theorems and surjectivity of invariant
differential operators on symmetric spaces and Lie groups, 
Bull.\ Amer.\ Math.\ Soc.\  {\bf 79} (1973), 129-132.

%\Lable{\rHeGG} {\it S.\ Helgason,}
%Groups and geometric analysis.
%American Mathematical Society 2000.

\Lable{\rHorm} {\it L.\ H\"ormander,}
The analysis of linear partial differential operators I.
Springer 1983.

\Lable{\rSchaefer} {\it H.\ H.\ Schaefer,}
Topological vector spaces. 2nd ed, Springer Verlag 1999. 

\Lable{\rSbook} {\it H.\ Schlichtkrull,}
Hyperfunctions and harmonic analysis on symmetric spaces,
Birk\-h\"au\-ser 1984.

\Lable{\rSsurvey}
{\it H.\ Schlichtkrull,}
The Paley-Wiener theorem for a reductive symmetric space, 
pp.\ 99-134 in
{\it Lie theory, Harmonic analysis on symmetric spaces -
General Plancherel theorems}, 
J.-P.\ Anker and B.\ \O rsted (eds.), Birkh\"auser Boston, 2005.

\bigskip\bigskip
\hbox{\hbox to 7 cm{\vbox{\hsize=7cm
E.P. van den Ban

Mathematisch Instituut

Universiteit Utrecht

P.O. Box 80010

3508 TA Utrecht

The Netherlands

E-mail: ban@math.uu.nl}}
\hbox to 6 cm{\vbox{\hsize=6cm
H. Schlichtkrull

Matematisk Institut

K\o benhavns Universitet

Universitetsparken 5

2100 K\o benhavn \O

Denmark

E-mail: schlicht@math.ku.dk}}}

\end